\documentclass[11pt]{article}

\usepackage{amsfonts,amssymb,graphics,epsfig,verbatim,bm,latexsym,amsmath,url,amsbsy,
geometry}

\begin{document}
  
\begin{center}

$\;$

\textbf{\Large Quasi Maximum-Likelihood Estimation \\
of Dynamic Panel Data Models}

\bigskip
\bigskip

Robert F. Phillips
\bigskip

Department of Economics \\
George Washington University \\
Washington, D.C. 20052\\
\bigskip
E-mail: \texttt{rphil@gwu.edu}
\bigskip

January 2017 \\

\bigskip
\bigskip

\end{center}

  \begin{abstract}
    This paper establishes the almost sure convergence and asymptotic normality of levels and 
differenced quasi maximum-likelihood (QML) estimators of dynamic panel data models. The QML estimators 
are robust with respect to initial conditions, 
conditional and time-series heteroskedasticity, and misspecification of the log-likelihood. 
The paper also provides an ECME algorithm for calculating  levels QML estimates. Finally, it uses Monte Carlo
experiments to compare the finite sample performance of levels and differenced QML estimators, 
the differenced GMM estimator, and the system GMM estimator. In these
experiments the QML estimators usually have smaller --- typically substantially smaller --- bias and root mean squared errors than the panel data
GMM estimators. 

  \end{abstract}
  
  \section{Introduction}
  
  Two prominent approaches to estimating a dynamic panel data model are generalized method of moments (GMM) and maximum likelihood (ML). Several authors have studied ML estimation of dynamic  panel data models; see, for example, Alvarez and Arellano (2004),
Anderson and Hsiao (1981), Hsiao et al. (2002), and Moral-Benito (2013),
among others.  As is well-known, the consistency and asymptotic normality of a ML estimator follows from ML theory assuming the likelihood is correctly specified
and standard regularity conditions are met.  On the other hand, strong distributional assumptions are not required to establish the sampling behavior of a GMM estimator. This fact would appear to make GMM more attractive than ML, but GMM has its drawbacks as well --- for example, GMM estimators are known to often have severe finite sample bias.  Furthermore, some papers have shown that the maximizer of a log-likelihood for a panel data model can be consistent and asymptotically normal under assumptions that do not require normality.  
Binder et al. (2005), for example, considered quasi-ML (QML) estimation of vector panel
autoregressions. Kruiniger (2013), on the other hand, studied QML estimation
of a first-order autoregressive (AR(1)) panel data model. And Phillips
(2010, 2015) examined QML estimation of a $p$th-order dynamic panel data
model. These papers provide conditions under which the log-likelihood for a
dynamic panel data model can be misspecified, and the maximizer of the quasi
log-likelihood is nevertheless consistent and asymptotically normal.

This paper makes several contributions to the literature on QML estimation. Like Phillips (2010, 2015), the model studied in this paper
includes $p$ lags of the dependent variable as well as other explanatory
variables. Phillips (2010, 1015), however, focused on QML estimation without
differencing the observations --- i.e., levels QML --- while assuming the
errors are unconditionally homoskedastic. The assumption of unconditional
homoskedasticity is more general than it might first appear, for it allows
for conditional heteroskedasticity. But it does not allow for time-series
heteroskedasticity. Allowing for more general forms of heteroskedasticity is
important, for QML estimation, although robust with respect to initial
conditions and misspecification of the log-likelihood, is not robust to
misspecification of the unconditional error variance-covariance matrix; see
also Alvarez and Arellano (2004). This paper, therefore, provides large $N$,
fixed $T$ asymptotics under more general conditions than those considered in
Phillips (2010, 2015) --- conditions that allow for time-series
heteroskedasticity. Indeed, the error variance-covariance matrix can be of a
general form.

Phillips (2010) provided a
straightforward iterative feasible generalized least squares algorithm for
calculating QML estimates when the errors in the dynamic regression model have an error-components structure. However, that procedure is not easily extended to
the case where the idiosyncratic errors are time-series heteroskedastic.
Furthermore, derivative-based algorithms can produce negative fitted variance
components when applied to error-components models if they are not substantially modified to avoid that outcome
(see also Meng and van Dyk 1998). This paper improves on these algorithms by providing an expectation conditional maximization either
(ECME) algorithm for calculating levels QML estimates that allows for
conditional and time-series heteroskedasticity. The ECME algorithm is
straightforward and guarantees non-negative estimated variance components.%

The paper also examines QML estimation after differencing the observations
(differenced QML). It shows that the ML estimator examined by
Hsiao et al. (2002) is consistent and asymptotically normal under more
general conditions than the conditions considered by Hsiao et al. (2002).
For example, Hsiao et al. (2002) assumed normality. This paper shows the estimator can be consistent and
asymptotically normal even if the log-likelihood is misspecified. Moreover,
restrictive initial conditions are not required, and the errors can be
conditionally heteroskedastic.

Finally, using simulated data, the finite sample behavior of levels and
differenced QML estimators are compared, and their finite sample behavior is
compared to the differenced GMM (Arellano and Bond 1991) and the system GMM
estimators (Blundell and Bond 1998). The Monte Carlo results show that,
compared to GMM estimators, the QML estimators have negligible finite sample
bias, and consequently they have smaller --- sometimes much smaller --- root
mean squared errors.
  
\section{QML via Regression Augmentation}\label{levels_qml}

Since Anderson and Hsiao (1981) it has been known that whether or not
application of ML estimation to a dynamic panel data model will yield a
consistent estimator as $N\rightarrow \infty $, with $T$ fixed, depends on
initial conditions. However, Phillips (2010) showed that, when QML
estimation is based on observations in levels (henceforth levels QML), it
does not depend on initial condition restrictions if the regression is
augmented with a suitable control function. This section extends the results
in Phillips (2010) by establishing the almost sure convergence and
asymptotic normality of levels QML estimation under weaker conditions than
thosed used in Phillips (2010). For example, the results provided here allow
for more general specifications of the error variance-covariance matrix.
This generalization is important because QML estimation is inconsistent if
the error variance-covariance matrix is misspecified.%

The model examined in this paper is the $p$th-order dynamic panel data model%
\begin{equation}
\boldsymbol{y}_{i}=\boldsymbol{Y}_{i}\delta _{0}+\boldsymbol{X}_{i}%
\boldsymbol{\beta }_{0}+\boldsymbol{e}_{i}\text{\ \ \ \ \ \ \ \ \ \ \ \ \ \
\ }\left( i=1,\ldots ,N\right) .  \label{model}
\end{equation}%
In this expression $\boldsymbol{y}_{i}=\left( y_{i1},\ldots ,y_{iT}\right)
^{\prime }$, $\boldsymbol{Y}_{i}=\left( \boldsymbol{y}_{i,-1},\ldots ,%
\boldsymbol{y}_{i,-p}\right) $, $\boldsymbol{y}_{i,-j}=\left(
y_{i,1-j},\ldots ,y_{i,T-j}\right) ^{\prime }$ ($j=1,\ldots ,p$), and $%
\boldsymbol{X}_{i}=\left( \boldsymbol{x}_{i1},\ldots ,\boldsymbol{x}%
_{iT}\right) ^{\prime }$, with $\boldsymbol{x}_{it}$ a $K\times 1$ vector of
explanatory variables that vary with $t$ (for at least some $i$). Moreover, $%
\boldsymbol{e}_{i}=\left( e_{i1},\ldots ,e_{iT}\right) ^{\prime }$ is a
vector of regression errors. For notational convenience, the numbering of
observed variables begins with $t=-p+1$.%

Straightforward ML estimation of the model in (\ref{model}) will not
generally yield a consistent estimator. To see why, let $\boldsymbol{y}%
_{i}^{o}=\left( y_{i0},\ldots ,y_{i,-p+1}\right) ^{\prime }$; let $%
\boldsymbol{x}_{i}$ be a column vector consisting of all of the distinct
elements of $\boldsymbol{x}_{i1},\ldots ,\boldsymbol{x}_{iT}$; and set $%
\boldsymbol{z}_{i}=\left( \boldsymbol{x}_{i}^{\prime },\boldsymbol{y}%
_{i}^{o\prime }\right) ^{\prime }$. Then, assuming $\boldsymbol{e}_{i}|%
\boldsymbol{z}_{i}\sim II\mathcal{N}\left( \mathbf{0},\Omega _{0}^{\ast
}\right) $, the log-likelihood is given by%
\begin{equation}
-\frac{NT}{2}\ln \left( 2\pi \right) -\frac{N}{2}\ln |\Omega ^{\ast }|-\frac{%
1}{2}\sum_{i=1}^{N}\boldsymbol{e}_{i}\left( \boldsymbol{\varphi }\right)
^{\prime }\Omega ^{\ast -1}\boldsymbol{e}_{i}\left( \boldsymbol{\varphi }%
\right) ,  \label{wrong_log-like}
\end{equation}%
where $\boldsymbol{e}_{i}\left( \boldsymbol{\varphi }\right) =\boldsymbol{y}%
_{i}-\boldsymbol{Y}_{i}\delta -\boldsymbol{X}_{i}\boldsymbol{\beta }$, and $%
\boldsymbol{\varphi }=\left( \boldsymbol{\delta }^{\prime },\boldsymbol{%
\beta }^{\prime }\right) ^{\prime }$. If $\Omega _{0}^{\ast }$ were known,
then maximizing the log-likelihood in (\ref{wrong_log-like}) yields the
generalized least squares (GLS)\ estimator based on $\Omega _{0}^{\ast }$,
and the consistency of that estimator requires $E\left( \boldsymbol{X}%
_{i}^{\prime }\Omega _{0}^{\ast -1}\boldsymbol{e}_{i}\right) =\mathbf{0}$
and $E\left( \boldsymbol{y}_{i,-j}^{\prime }\Omega _{0}^{\ast -1}\boldsymbol{%
e}_{i}\right) =\mathbf{0}$ ($j=1,\ldots ,p$).%

We have $E\left( \boldsymbol{X}_{i}^{\prime }\Omega _{0}^{\ast -1}%
\boldsymbol{e}_{i}\right) =\mathbf{0}$ if the regressors in $\boldsymbol{X}%
_{i}$ are strictly exogenous with respect to the errors in $\boldsymbol{e}%
_{i}$. But the moment restrictions $E\left( \boldsymbol{y}_{i,-j}^{\prime
}\Omega _{0}^{\ast -1}\boldsymbol{e}_{i}\right) =\mathbf{0}$ ($j=1,\ldots ,p$%
) depends on an even stronger assumption, which is summarized in Lemma 1.%

\bigskip

\noindent \textbf{Lemma 1}. If $E\left( \boldsymbol{e}_{i}\boldsymbol{y}%
_{i}^{o\prime }\right) =\mathbf{0}$, $E\left( \boldsymbol{e}_{i}\boldsymbol{x%
}_{i}^{\prime }\right) =\mathbf{0}$, and $E\left( \boldsymbol{e}_{i}%
\boldsymbol{e}_{i}^{\prime }\right) =\Omega _{0}^{\ast }$, then $E\left( 
\boldsymbol{y}_{i,-j}^{\prime }\Omega _{0}^{\ast -1}\boldsymbol{e}%
_{i}\right) =\mathbf{0}$ ($j=1,\ldots ,p$).%

\bigskip

\noindent \textbf{Proof}. See Appendix A.

\bigskip

According to Lemma 1, if the regressors in $\boldsymbol{x}_{it}$ 
\emph{and} the initial values of the dependent variable $y_{i0},\ldots
,y_{i,-p+1}$ are uncorrelated with the errors $e_{i1},\ldots ,e_{iT}$, then $%
E\left( \boldsymbol{y}_{i,-j}^{\prime }\Omega _{0}^{\ast -1}\boldsymbol{e}%
_{i}\right) =\mathbf{0}$ ($j=1,\ldots ,p$). However, assuming the initial
values of the dependent variable are uncorrelated with subsequent errors is
quite restrictive. For example, a commonly used model for the errors is the
error-components model%
\begin{equation}
e_{it}=c_{i}+v_{it}.  \label{EC}
\end{equation}%
If the $v_{it}$s are uncorrelated, we can take $v_{it}$ to be uncorrelated
with the elements of $\boldsymbol{y}_{i}^{o}$, for $t\geq 1$, but assuming
the elements of $\boldsymbol{y}_{i}^{o}$ are also uncorrelated with $c_{i}$
is a strong initial condition restriction.%

Fortunately, we need make no such initial condition assumption if the model
in (\ref{model}) is augmented with a suitable control function. Nor need we
assume the regressors in $\boldsymbol{x}_{it}$ are strictly exogenous with
respect to the $e_{it}$s. The possible correlation between the elements in $%
\boldsymbol{e}_{i}$ and the elements in $\boldsymbol{z}_{i}$ can be
controlled for by the linear projection of $e_{it}$ on 1 and $\boldsymbol{z}%
_{i}$:%
\begin{equation}
e_{it}=\mu _{0}+\boldsymbol{z}_{i}^{\prime }\boldsymbol{\theta }_{0}+u_{it},%
\text{ \ \ \ \ \ }(t=1,\ldots ,T\text{, \ }i=1,\ldots ,N)  \label{lp}
\end{equation}%
where $\boldsymbol{\theta }_{0}=Var\left( \boldsymbol{z}_{i}\right)^{-1}Cov\left( \boldsymbol{z}_{i},e_{it}\right)$ and
$\mu _{0}=E\left( e_{it}\right)-E\left( \boldsymbol{z}_{i}\right)^{\prime }\boldsymbol{\theta }_{0}$.

The linear projection parameters $\mu _{0}$ and $\boldsymbol{\theta }%
_{0}$ exist and depend on neither $i$ nor $t$ if $E\left( e_{it}\right)$ and the moments in $Cov\left( \boldsymbol{z}_{i},e_{it}\right)$ depend on neither $i$ nor $t$ and the moments in $Var\left( \boldsymbol{z}_{i}\right)$ and  
$E\left( \boldsymbol{z}_{i}\right)$ do not depend on $i$. The
restriction that the linear projection parameters are independent of $t$ is
met if the errors have a one-way error-components structure given by (\ref%
{EC}) and $v_{it}$ is a mean zero random variable that is uncorrelated with
the elements of $\boldsymbol{z}_{i}$ for $t\geq 1$. Then $Cov\left( 
\boldsymbol{z}_{i},e_{it}\right) =Cov\left( \boldsymbol{z}_{i},c_{i}\right) $
and $E\left( e_{it}\right) =E\left( c_{i}\right) $ for $t\geq 1$. For this
case, the linear projection reduces to that considered in Phillips (2010,
2015). Specifically, we have 
\begin{equation}
c_{i}=\mu _{0}+\boldsymbol{z}_{i}^{\prime }\boldsymbol{\theta }_{0}+a_{i}%
\text{ \ \ \ \ \ }(i=1,\ldots ,N)  \label{lp_ec}
\end{equation}%
(cf Phillips 2010, p. 411, Eq. (2)).\footnote{%
See also Chamberlain (1982, 1984) and Kruiniger (2013), who uses a linear
projection of an individual effect on $y_{i0}$. The linear projection
parameters used in Kruiniger (2013) are implicitly assumed to be independent
of $i$.} If the errors can be decomposed as in Eq. (\ref{EC}), then $\mu _{0}+%
\boldsymbol{z}_{i}^{\prime }\boldsymbol{\theta }_{0}$ controls for possible
correlation between time-invariant unobservables, captured by $c_{i}$, and
the elements of $\boldsymbol{z}_{i}$. 

Another, albeit trivial, case in which the linear projection parameters
depend on neither $i$ nor $t$ is when there are no individual specific
effects and the $e_{it}$s are uncorrelated among themselves and with the
elements of $\boldsymbol{z}_{i}$, for $t\geq 1$. In this case, $\boldsymbol{%
\theta }_{0}=\mathbf{0}$, and the linear projection in (\ref{lp}) simplifies
to $e_{it}=\mu _{e}+u_{it}$, where $E\left( e_{it}\right) =\mu _{e}$. This
example illustrates that the necessity of adding the control function $\mu
_{0}+\boldsymbol{z}_{i}^{\prime }\boldsymbol{\theta }_{0}$ follows from the
presence of unobservable time-invariant omitted variables, which are
captured by $c_{i}$.

Moreover, although it is obvious we must include $\boldsymbol{x}_{i}$ in the
control function when the regressors in $\boldsymbol{x}_{it}$ are correlated
with $c_{i}$, it is also true that we typically must do so even when all of
the regressors in $\boldsymbol{x}_{it}$ are uncorrelated $c_{i}$, as in the
random effects model. To see this, consider the linear projection of $c_{i}$
on just 1 and $\boldsymbol{y}_{i}^{o}$:%
\begin{equation}
c_{i}=\mu _{y0}+\boldsymbol{y}_{i}^{o\prime }\boldsymbol{\theta }_{y0}+a_{yi}%
\text{ \ \ \ \ \ }(i=1,\ldots ,N),  \label{wrong_lp}
\end{equation}%
where $\boldsymbol{\theta }_{y0}=Var\left( \boldsymbol{y}_{i}^{o}\right)
^{-1}Cov\left( \boldsymbol{y}_{i}^{o},c_{i}\right) $ and $\mu _{y0}=E\left(
c_{i}\right) -E\left( \boldsymbol{y}_{i}^{o^{\prime }}\right) \boldsymbol{%
\theta }_{y0}$. If we augment the model in (\ref{model}) with the control
function $\mu _{y0}+\boldsymbol{y}_{i}^{o^{\prime }}\boldsymbol{\theta }%
_{y0} $ rather than the control function $\mu _{0}+\boldsymbol{z}%
_{i}^{\prime }\boldsymbol{\theta }_{0}$, then the error term in the
augmented model is $a_{yi}+v_{it}$ rather than $a_{i}+v_{it}$, and, in order
for QML estimation of the augmented model to be consistent, we must have not
just $Cov\left( \boldsymbol{y}_{i}^{o},a_{yi}\right) =\mathbf{0}$, which the
linear projection in (\ref{wrong_lp}) ensures, but also $Cov\left( 
\boldsymbol{x}_{i},a_{yi}\right) =\mathbf{0}$, which the linear projection
in (\ref{wrong_lp}) does not guarantee. Indeed, given $Cov\left( \boldsymbol{%
x}_{i},c_{i}\right) =\mathbf{0}$, the result $Cov\left( \boldsymbol{x}%
_{i},a_{yi}\right) =\mathbf{0}$ is not guaranteed unless $Cov\left( 
\boldsymbol{x}_{i},\boldsymbol{y}_{i}^{o^{\prime }}\boldsymbol{\theta }%
_{y0}\right) =\mathbf{0}$,\footnote{%
This conclusion follows from $Cov\left( \boldsymbol{x}_{i},a_{yi}\right)
=Cov(\boldsymbol{x}_{i},c_{i}-\mu _{y0}-\boldsymbol{y}_{i}^{o\prime }%
\boldsymbol{\theta }_{y0})=-Cov(\boldsymbol{x}_{i},\boldsymbol{y}%
_{i}^{o\prime }\boldsymbol{\theta }_{y0})$ if $Cov\left( \boldsymbol{x}%
_{i},c_{i}\right) =\mathbf{0}$.} which will not be satisfied in general
assuming $\boldsymbol{\theta }_{y0}\neq \mathbf{0}$.%

This last example illustrates that results obtained for the AR(1) panel data
model (see Kruiniger 2013) or the AR(p) panel data model (see Alvarez and
Arellano 2004) do not extend in a straightforward manner to models with
additional regressors even under the random effects assumption that the
elements of $\boldsymbol{x}_{it}$ are uncorrelated with $c_{i}$. For
example, in his treatment of the \textquotedblleft random
effects\textquotedblright\ case of the AR(1) panel data model, Kruiniger
includes a linear projection of $c_{i}$ on the initial value $y_{i0}$ in a
control function. However, such a control function will not suffice if there
are additional regressors even when these additional regressors are
uncorrelated with $c_{i}$.

Equations (\ref{model}) and (\ref{lp}) imply the augmented dynamic
panel data model%
\begin{equation}
\boldsymbol{y}_{i}=\boldsymbol{W}_{i}\boldsymbol{\gamma }_{0}+\boldsymbol{u}%
_{i},\text{\ \ \ \ \ \ \ \ \ \ \ \ \ \ \ }\left( i=1,\ldots ,N\right) ,
\label{aug_model}
\end{equation}%
where $\boldsymbol{W}_{i}=\left( \boldsymbol{Y}_{i},\boldsymbol{Z}%
_{i}\right) $, $\boldsymbol{Z}_{i}=\left( \boldsymbol{X}_{i},\boldsymbol{%
\iota }\,,\boldsymbol{\iota \,z}_{i}^{\prime }\right) $, $\boldsymbol{\iota }
$ is a $T\times 1$ vector of ones, and $\boldsymbol{\gamma }_{0}=\left( 
\boldsymbol{\delta }_{0}^{\prime },\boldsymbol{\beta }_{0}^{\prime },\mu
_{0},\boldsymbol{\theta }_{0}^{\prime }\right) ^{\prime }$. The errors in
this augmented model --- $\boldsymbol{u}_{i}=\left( u_{i1},\ldots
,u_{iT}\right) ^{\prime }$ --- are now uncorrelated with the elements of $%
\boldsymbol{Z}_{i}$ by construction. Thus, upon letting $\Omega _{0}=E\left( 
\boldsymbol{u}_{i}\boldsymbol{u}_{i}^{\prime }\right) $, we have $E\left( 
\boldsymbol{Z}_{i}^{\prime }\Omega _{0}^{-1}\boldsymbol{u}_{i}\right) =%
\mathbf{0}$. Moreover, because $E\left( \boldsymbol{u}_{i}\boldsymbol{y}%
_{i}^{o\prime }\right) =\mathbf{0}$ and $E\left( \boldsymbol{u}_{i}%
\boldsymbol{x}_{i}^{\prime }\right) =\mathbf{0}$, it follows from Lemma 1
that $E\left( \boldsymbol{y}_{i,-j}^{\prime }\Omega _{0}^{-1}\boldsymbol{u}%
_{i}\right) =\mathbf{0}$ ($j=1,\ldots ,p$). The preceding shows $E\left( 
\boldsymbol{W}_{i}^{\prime }\Omega _{0}^{-1}\boldsymbol{u}_{i}\right) =%
\mathbf{0}$.

Now consider the quasi log-likelihood for the augmented model in (\ref%
{aug_model}): $\sum_{i=1}^{N}l_{i}\left( \boldsymbol{\psi }\right) $, where 
\begin{equation*}
l_{i}\left( \boldsymbol{\psi }\right) =-\frac{T}{2}\ln \left( 2\pi \right) -%
\frac{1}{2}\ln |\Omega |-\frac{1}{2}\boldsymbol{u}_{i}\left( \boldsymbol{%
\gamma }\right) ^{\prime }\Omega ^{-1}\boldsymbol{u}_{i}\left( \boldsymbol{%
\gamma }\right) ,
\end{equation*}%
$\boldsymbol{u}_{i}\left( \boldsymbol{\gamma }\right) =\boldsymbol{y}_{i}-%
\boldsymbol{W}_{i}\boldsymbol{\gamma }$, $\boldsymbol{\gamma }=\left( 
\boldsymbol{\delta }^{\prime },\boldsymbol{\beta }^{\prime },\mu ,%
\boldsymbol{\theta }^{\prime }\right) ^{\prime }$, $\boldsymbol{\psi }%
=\left( \boldsymbol{\gamma }^{\prime },\boldsymbol{\omega }^{\prime }\right)
^{\prime }$, $\boldsymbol{\omega }=$ vech$\left( \Omega \right) $, and $%
\Omega $ is a positive definite matrix. For known $\boldsymbol{\omega }_{0}=$
vech$\left( \Omega _{0}\right) $, the maximizer of this log-likelihood is
the GLS estimator $\widehat{\boldsymbol{\gamma }}_{GLS}=\left( \sum_{i=1}^{N}%
\boldsymbol{W}_{i}^{\prime }\Omega _{0}^{-1}\boldsymbol{W}_{i}\right)
^{-1}\sum_{i=1}^{N}\boldsymbol{W}_{i}^{\prime }\Omega _{0}^{-1}\boldsymbol{y}%
_{i}$, and this estimator is consistent because $E\left( \boldsymbol{W}%
_{i}^{\prime }\Omega _{0}^{-1}\boldsymbol{u}_{i}\right) =\mathbf{0}$.
Moreover, if $\widehat{\Omega }$ is a consistent estimator of $\Omega _{0},$
the feasible GLS (FGLS) estimator $\widehat{\boldsymbol{\gamma }}%
_{FGLS}=\left( \sum_{i=1}^{N}\boldsymbol{W}_{i}^{\prime }\widehat{\Omega }%
^{-1}\boldsymbol{W}_{i}\right) ^{-1}\sum_{i=1}^{N}\boldsymbol{W}_{i}^{\prime
}\widehat{\Omega }^{-1}\boldsymbol{y}_{i}$ is also consistent.%

However, the large $N$ (fixed $T$) distribution of such a FGLS estimator
depends on the first-round estimator of $\boldsymbol{\gamma }_{0}$ used to
estimate $\Omega _{0}$ (see Phillips 2010). An alternative that does not
depend on a first-round estimator is to estimate $\boldsymbol{\psi }%
_{0}=\left( \boldsymbol{\gamma }_{0}^{\prime },\boldsymbol{\omega }%
_{0}^{\prime }\right) ^{\prime }$ by maximizing the quasi log-likelihood $%
\sum_{i=1}^{N}l_{i}\left( \boldsymbol{\psi }\right) $.%

Theorems 1 and 2 provide sufficient conditions for the almost sure
convergence of the QML estimator and its asymptotic normality (as $%
N\rightarrow \infty $, with $T$ fixed). In order to state the theorems, set $%
L_{N}\left( \boldsymbol{\psi }\right) =N^{-1}\sum_{i=1}^{N}l_{i}\left( 
\boldsymbol{\psi }\right) $ and $\boldsymbol{H}_{N}\left( \boldsymbol{\psi }%
\right) =\partial ^{2}L_{N}\left( \boldsymbol{\psi }\right) /\partial 
\boldsymbol{\psi }\partial \boldsymbol{\psi }^{\prime }$; let $x_{itk}$
denote the $k$th element of $\boldsymbol{x}_{it}$; and set $\Psi =\left\{ 
\boldsymbol{\psi }=\left( \boldsymbol{\gamma }^{\prime },\boldsymbol{\omega }%
^{\prime }\right) ^{\prime }\boldsymbol{\in 
\mathbb{R}
}^{m}:\Omega \text{ is positive definite}\right\} $.%

\bigskip

\sloppy
\noindent \textbf{Theorem 1}. Assume the following conditions are satisfied:

\begin{description}
\item[C1:] $E\left\vert y_{it}\right\vert ^{2+\epsilon }<M$ and $E\left\vert
x_{itk}\right\vert ^{2+\epsilon }<M$ for all $i$, $t$, and $k$ and some $%
\epsilon >0$ and $M<\infty $;

\item[C2:] $Var\left( \boldsymbol{z}_{i}\right) =\Xi _{zz}$ for all $i$,
with $\Xi _{zz}$ a positive definite matrix, $E\left( \boldsymbol{z}%
_{i}\right) =\boldsymbol{\mu }_{z}$ for all $i$, and $E\left( e_{it}\right)
=\mu _{e}$ and $E\left( \boldsymbol{z}_{i}e_{it}\right) =\boldsymbol{\varrho 
}_{ze}$ for all $i$ and $t\geq 1$;

\item[C3:] $E\left( \boldsymbol{u}_{i}\boldsymbol{u}_{i}^{\prime }\right)
=\Omega _{0}$ for all $i$, with $\Omega _{0}$ a positive definite matrix;

\item[C4:] the limits $\lim_{N\rightarrow \infty }N^{-1}\sum_{i}E\left(
y_{is}y_{it}\right) $, $\lim_{N\rightarrow \infty }N^{-1}\sum_{i}E\left(
y_{is}x_{itk}\right) $, and $\lim_{N\rightarrow \infty
}N^{-1}\sum_{i}E\left( x_{isj}x_{itk}\right) $ exist for all $s$, $t$, $j$,
and $k$; and

\item[C5:] the vectors $\left( \boldsymbol{z}_{1}^{\prime },\boldsymbol{y}%
_{1}^{\prime }\right) ^{\prime },\ldots ,\left( \boldsymbol{z}_{N}^{\prime },%
\boldsymbol{y}_{N}^{\prime }\right) ^{\prime }$ are independent for all $N$.
\end{description}

\fussy
\noindent Then $E\left[\partial L_{N}\left( \boldsymbol{\psi }_{0}\right) /\partial 
\boldsymbol{\psi }\right] = \bold{0}$ and the limit $\boldsymbol{H}\left( \boldsymbol{\psi }\right)
=\lim_{N\rightarrow \infty }E\left[ \boldsymbol{H}_{N}\left( \boldsymbol{%
\psi }\right) \right] $ exists. Moreover, if $\boldsymbol{H}_{0}=\boldsymbol{%
H}\left( \boldsymbol{\psi }_{0}\right) $ is negative definite, then there is
a compact subset, say $\overline{\Psi }$, of $\Psi $, with $\boldsymbol{\psi 
}_{0}$ in its interior, and there is a measurable maximizer, $\widehat{%
\boldsymbol{\psi }}$, of $L_{N}\left( \cdot \right) $ in $\overline{\Psi }$
such that $\widehat{\boldsymbol{\psi }}\overset{a.s.}{\rightarrow }%
\boldsymbol{\psi }_{0}$ ($N\rightarrow \infty $, $T$ fixed).%

\bigskip

\noindent \textbf{Proof}. See Appendix B.%

\bigskip

\noindent \textbf{Theorem 2}. Assume Conditions C2--C5 are satisfied, $%
\boldsymbol{H}_{0}$ is negative definite, and the following conditions are
satisfied:

\begin{description}
\item[C1$^{\prime }$:] $E\left\vert y_{it}\right\vert ^{4+\epsilon }<M$ and $%
E\left\vert x_{itk}\right\vert ^{4+\epsilon }<M$ for all $i$, $t$, and $k$
and some $\epsilon >0$ and $M<\infty $; and

\item[C6:] the limit $\mathcal{I}_{0}=\lim_{N\rightarrow \infty
}N^{-1}\sum_{i}E\left[ \left( \partial l_{i}\left( \boldsymbol{\psi }%
_{0}\right) /\partial \boldsymbol{\psi }\right) \left( \partial l_{i}\left( 
\boldsymbol{\psi }_{0}\right) /\partial \boldsymbol{\psi }\right) ^{\prime }%
\right] $ exists and is positive definite.
\end{description}

\noindent Then $\sqrt{N}\left( \widehat{\boldsymbol{\psi }}-\boldsymbol{\psi 
}_{0}\right) \overset{d}{\rightarrow }\mathcal{N}\left( \mathbf{0},%
\boldsymbol{H}_{0}^{-1}\mathcal{I}_{0}\boldsymbol{H}_{0}^{-1}\right) $ ($%
N\rightarrow \infty $, $T$ fixed).%

\bigskip

\noindent \textbf{Proof}. See Appendix C.%

\bigskip

In order for the QML estimator to be consistent and asymptotically normal, it must be the case that
the true parameter vector, $ \boldsymbol{\psi }_{0}$, uniquely maximizes the expected log-likelihood, at least within a neighborhood of $ \boldsymbol{\psi }_{0}$.  Conditions C1 through C3 are mild, and they suffice to guarantee that $ \boldsymbol{\psi }_{0}$ is indeed a stationary value of the expected log-likelihood.  But the fact that $\boldsymbol{\psi }_{0}$ is a stationary value is necessary but not sufficient to ensure it is a unique maximizer of the expected log-likelihood. The matrix $\boldsymbol{H}_{0}$ must also be negative definite. If the log-likelihood $\sum_{i=1}^{N}l_{i}%
\left( \boldsymbol{\psi }\right) $ is correctly specified, that is, if $%
\boldsymbol{u}_{i}$ is normally distributed with mean vector $\mathbf{0}$
and variance-covariance matrix $\Omega _{0}$, conditionally on $\boldsymbol{z%
}_{i}$, then by well-known ML theory, we have $\boldsymbol{H}_{0}=-\mathcal{I%
}_{0}$, and $\boldsymbol{H}_{0}$ exists and is negative definite
by virtue of Condition C6. However, even when $\sum_{i=1}^{N}l_{i}\left( 
\boldsymbol{\psi }\right) $ is misspecified,  $\boldsymbol{H}_{0}$ can be
shown to be negative definite in particular cases. Phillips (2015), for example, provides an example in which $\boldsymbol{H}_{0}$ is negative definite under conditions that do not include
normality.

Moreover, $\Omega _{0}$ is the unconditional variance-covariance matrix of $%
\boldsymbol{u}_{i}$, and, although it does not depend on $i$, the
variance-covariance matrix of $\boldsymbol{u}_{i}$ conditionally on $%
\boldsymbol{z}_{i}$ may depend on $i$ --- for example, the errors may be
conditionally heteroskedastic (see also Phillips 2010, 2015). The errors can
also be unconditionally time-series heteroskedastic, for the diagonal elements of $\Omega _{0}$ can differ.%

Furthermore, the conditions in Theorems 1 and 2 do not require the
random vectors $\left( \boldsymbol{z}_{1}^{\prime },\boldsymbol{y}%
_{1}^{\prime }\right) ^{\prime },\ldots ,\left( \boldsymbol{z}_{N}^{\prime },%
\boldsymbol{y}_{N}^{\prime }\right) ^{\prime }$ be drawn from a common
distribution. On the other hand, Conditions C2 and C3 imply some homogeneity is required. 

Estimators previously considered in the literature are covered by
Theorems 1 and 2. Blundell and Bond (1998) considered a conditional GLS
estimator of an AR(1) panel data model that relied on augmenting the
regression model with the initial observation on the dependent variable.
They argued that if the error components are homoskedastic across
individuals and time, then restrictions on the initial conditions can be
used to derive the GLS estimator. Theorems 1 and 2, however, show that these
conditions are unnecessarily restrictive. The errors can be conditionally
and time-series heteroskedastic. Moreover, initial condition restrictions
are not needed. All that is required is that the moments defining the
control function parameters exist and depend on on neither $i$ nor $t$.
Furthermore, the structured error variance-covariance matrices, such as
those considered by Phillips (2010, 2015) and Kruiniger (2013), are special
cases of $\Omega _{0}$, and, therefore, Theorems 1 and 2 cover those
cases.

 \section{Fixed-Effects QML}\label{differenced_qml}

An alternative to first augmenting the regression model with a control
function and then applying QML estimation to the model in levels is to
instead first difference the observations and then apply QML estimation. In
the literature, ML or QML estimation based on first differencing the
observations has been referred to as fixed-effects ML estimation (e.g.,
Hsiao et al. 2002) or fixed-effects QML estimation (e.g., Kruiniger 2013).
This description, however, should not lead one to interpret levels QML
estimation as random-effects QML, for the results in Section \ref{levels_qml}
make clear that levels QML estimation is not restricted to random-effects
models with regressors that are exogenous with respect to $c_{i}$.

Kruiniger (2013) studied differenced QML for an AR(1) panel data model.
Hsiao et al. (2002), on the other hand, studied ML\ estimation, after
differencing, and, like this paper, considered a model with additional
explanatory variables beyond a lagged dependent variable. This section shows
that likelihood-based methods using differences are consistent and
asymptotically normal under much weaker conditions than those assumed in
Hsiao et al. (2002).

\sloppy Instead of augmenting the regression with a control function that
involves $\boldsymbol{y}_{i}^{o}$, differenced QML requires estimation of a
system of equations that includes a separate linear projection for each
initial difference $\Delta y_{i,-p+2},\ldots ,\Delta y_{i1}$, where $\Delta
y_{it}=y_{it}-y_{i,t-1}$. Specifically, suppose $Var\left( \boldsymbol{x}%
_{i}\right) $ is positive definite, and set $\boldsymbol{\theta }%
_{0,p+1-j}=Var\left( \boldsymbol{x}_{i}\right) ^{-1}Cov\left( \boldsymbol{x}%
_{i},\Delta y_{i,-j+2}\right) $ and $\mu _{0,p+1-j}=E\left( \Delta
y_{i,-j+2}\right) -E\left( \boldsymbol{x}_{i}^{\prime }\right) \boldsymbol{%
\theta }_{0,p+1-j}$ ($j=1,\ldots ,p$). Then, system differenced QML relies
on estimating the linear projections%
\begin{equation}
\Delta y_{i,-j+2}=\mu _{0,p+1-j}+\boldsymbol{x}_{i}^{\prime }\boldsymbol{%
\theta }_{0,p+1-j}+r_{i,p+1-j}\text{ \ \ \ \ \ \ \ \ }(j=1,\ldots ,p).
\label{diff_lp}
\end{equation}%
Here $r_{i,p+1-j}$ is a linear projection residual, which is, by
construction, uncorrelated with all of the elements of $\boldsymbol{x}_{i}$.
Note that because the linear projection in (\ref{diff_lp}) does not specify
how $\Delta y_{i,-j+2}$ was generated it does not depend on initial
condition restrictions. In addition to the linear projection equations in (%
\ref{diff_lp})\ we also estimate the differenced equation:%
\begin{equation}
\Delta \boldsymbol{y}_{i}=\Delta \boldsymbol{Y}_{i}\boldsymbol{\delta }%
_{0}+\Delta \boldsymbol{X}_{i}\boldsymbol{\beta }_{0}+\Delta \boldsymbol{e}%
_{i}\text{\ \ \ \ \ \ \ \ }\left( i=1,\ldots ,N\right) ,  \label{diff_model}
\end{equation}%
where $\Delta \boldsymbol{y}_{i}=\left( \Delta y_{i2},\ldots ,\Delta
y_{iT}\right) ^{\prime }$, $\Delta \boldsymbol{Y}_{i}=\left( \Delta 
\boldsymbol{y}_{i,-1},\ldots ,\Delta \boldsymbol{y}_{i,-p}\right) $, and $%
\Delta \boldsymbol{y}_{i,-j}=\left( \Delta y_{i,-j+2},\ldots \Delta
y_{i,T-j}\right) ^{\prime }$ ($j=1,\ldots ,p$). Moreover, $\Delta 
\boldsymbol{X}_{i}=\left( \Delta \boldsymbol{x}_{i2},\ldots ,\Delta 
\boldsymbol{x}_{iT}\right) ^{\prime }$, $\Delta \boldsymbol{x}_{it}=%
\boldsymbol{x}_{it}-\boldsymbol{x}_{i,t-1}$, and $\Delta \boldsymbol{e}%
_{i}=\left( \Delta e_{i2},\ldots ,\Delta e_{iT}\right) ^{\prime }$, with $%
\Delta e_{it}=e_{it}-e_{i,t-1}$. For differenced QML, the equations in (\ref%
{diff_lp}) and (\ref{diff_model}) are estimated as a system given by%
\begin{equation}
\widetilde{\boldsymbol{y}}_{i}=\widetilde{\boldsymbol{W}}_{i}\boldsymbol{%
\eta }_{0}+\widetilde{\boldsymbol{u}}_{i}\text{\ \ \ \ \ \ \ \ \ \ \ \ \ \ \ 
}\left( i=1,\ldots ,N\right) ,  \label{levels_sys}
\end{equation}%
with $\widetilde{\boldsymbol{y}}_{i}=\left( \Delta y_{i,-p+2},\ldots ,\Delta
y_{i1},\Delta \boldsymbol{y}_{i}^{\prime }\right) ^{\prime }$, $\widetilde{%
\boldsymbol{u}}_{i}=\left( r_{i1},\ldots ,r_{ip},\Delta \boldsymbol{e}%
_{i}^{\prime }\right) ^{\prime },$%
\begin{equation*}
\widetilde{\boldsymbol{W}}_{i}=\left( 
\begin{array}{ccc}
\mathbf{0} & \mathbf{0} & \boldsymbol{I}_{p}\otimes \left( 1,\boldsymbol{x}%
_{i}^{\prime }\right) \\ 
\Delta \boldsymbol{Y}_{i} & \Delta \boldsymbol{X}_{i} & \mathbf{0}%
\end{array}%
\right) ,
\end{equation*}%
and $\boldsymbol{\eta }_{0}=\left( \boldsymbol{\delta }_{0}^{\prime },%
\boldsymbol{\beta }_{0}^{\prime },\mu _{01},\boldsymbol{\theta }%
_{01}^{\prime },\mu _{02},\boldsymbol{\theta }_{02}^{\prime },\ldots ,\mu
_{0p},\boldsymbol{\theta }_{0p}^{\prime }\right) ^{\prime }$.%

If $\widetilde{\boldsymbol{u}}_{i}$ is multivariate normal with mean vector $%
\mathbf{0}$ and variance-covariance matrix $\Upsilon _{0}$ conditional on $%
\boldsymbol{x}_{i}$, then the log-likelihood for the system in (\ref%
{levels_sys})\ is $\sum_{i=1}^{N}\widetilde{l}_{i}\left( \boldsymbol{\lambda 
}\right) $, where%
\begin{equation*}
\widetilde{l}_{i}\left( \boldsymbol{\lambda }\right) =-\frac{\left(
T+p-1\right) }{2}\ln \left( 2\pi \right) -\frac{1}{2}\ln |\Upsilon |-\frac{1%
}{2}\widetilde{\boldsymbol{u}}_{i}\left( \boldsymbol{\eta }\right) ^{\prime
}\Upsilon ^{-1}\widetilde{\boldsymbol{u}}_{i}\left( \boldsymbol{\eta }%
\right) ,
\end{equation*}%
$\widetilde{\boldsymbol{u}}_{i}\left( \boldsymbol{\eta }\right) =\widetilde{%
\boldsymbol{y}}_{i}-\widetilde{\boldsymbol{W}}_{i}\boldsymbol{\eta }$, $%
\boldsymbol{\eta }=\left( \boldsymbol{\delta }^{\prime },\boldsymbol{\beta }%
^{\prime },\mu _{1},\boldsymbol{\theta }_{1}^{\prime },\mu _{2},\boldsymbol{%
\theta }_{2}^{\prime },\ldots ,\mu _{p},\boldsymbol{\theta }_{p}^{\prime
}\right) ^{\prime }$, $\boldsymbol{\lambda }=\left( \boldsymbol{\eta }%
^{\prime },\boldsymbol{\upsilon }^{\prime }\right) ^{\prime }$, and $%
\boldsymbol{\upsilon }=$ vech$\left( \Upsilon \right) $. Also, set $%
\widetilde{L}_{N}\left( \boldsymbol{\lambda }\right) =N^{-1}\sum_{i=1}^{N}%
\widetilde{l}_{i}\left( \boldsymbol{\lambda }\right) $, $\widetilde{%
\boldsymbol{H}}_{N}\left( \boldsymbol{\lambda }\right) =\partial ^{2}%
\widetilde{L}_{N}\left( \boldsymbol{\lambda }\right) /\partial \boldsymbol{%
\lambda }\partial \boldsymbol{\lambda }^{\prime }$, and $\Lambda =\left\{ 
\boldsymbol{\lambda }=\boldsymbol{\left( \boldsymbol{\eta }^{\prime },%
\boldsymbol{\upsilon }^{\prime }\right) ^{\prime }\in 
\mathbb{R}
}^{n}:\Upsilon \text{ is positive definite}\right\} $.%

\fussy The maximizer of $\sum_{i=1}^{N}\widetilde{l}_{i}\left( \cdot \right) 
$ is a ML estimator given normality, but even if the log-likelihood is
misspecified --- that is, the errors are not normally distributed given $%
\boldsymbol{x}_{i}$, nor are they necessarily conditionally homoskedastic
--- maximizing $\sum_{i=1}^{N}\widetilde{l}_{i}\left( \cdot \right) $ will
still yield a consistent and asymptotically normal estimator under suitable
conditions. Sufficient conditions are provided in Theorems 3 and 4.%

\bigskip

\noindent \textbf{Theorem 3}. Suppose C1, C4, and C5 are satisfied. Further
assume:

\begin{description}
\item[C2$^{\prime }$:] $Var\left( \boldsymbol{x}_{i}\right) =\Xi _{xx}$ for
all $i$, with $\Xi _{xx}$ positive definite, $E\left( \boldsymbol{x}%
_{i}\right) =\boldsymbol{\mu }_{x}$ for all $i$, $E\left( \Delta
y_{i,-j+2}\right) =\mu _{\Delta y_{j}}$ and $E\left( \boldsymbol{x}%
_{i}\Delta y_{i,-j+2}\right) =\boldsymbol{\varrho }_{x\Delta y_{j}}$ for all 
$i$ $(j=1,\ldots ,p)$, and $Cov\left( \boldsymbol{x}_{i},\Delta \boldsymbol{e%
}_{i}\right) =\mathbf{0}$; also,

\item[C3$^{\prime }$:] $E\left( \widetilde{\boldsymbol{u}}_{i}\widetilde{%
\boldsymbol{u}}_{i}^{\prime }\right) =\Upsilon _{0}$ for all $i$, with $%
\Upsilon _{0}$ a positive definite matrix.
\end{description}

\noindent Then $E\left[\partial \widetilde{L}_{N}\left( \boldsymbol{\lambda }_{0}\right) /\partial 
\boldsymbol{\lambda }\right] = \bold{0}$, where $\boldsymbol{\lambda }_{0}=\left( 
\boldsymbol{\eta }_{0}^{\prime },\boldsymbol{\upsilon }_{0}^{\prime }\right)
^{\prime }$ and $\boldsymbol{\upsilon }_{0}=$ vech$\left( \Upsilon
_{0}\right) $.  Furthermore, the limit $\widetilde{\boldsymbol{H}}\left( \boldsymbol{%
\lambda }\right) =\lim_{N\rightarrow \infty }\widetilde{\boldsymbol{H}}%
_{N}\left( \boldsymbol{\lambda }\right) $ exists. Moreover, if $\widetilde{%
\boldsymbol{H}}_{0}=\widetilde{\boldsymbol{H}}\left( \boldsymbol{\lambda }%
_{0}\right) $ is negative definite, there is a compact subset, say $\overline{\Lambda }$, of $%
\Lambda $, with $\boldsymbol{\lambda }_{0}$ in its interior, and there is a
measurable maximizer, $\widehat{\boldsymbol{\lambda }}$, of $\widetilde{L}%
_{N}\left( \cdot \right) $ in $\overline{\Lambda }$ such that $\widehat{%
\boldsymbol{\lambda }}\overset{a.s.}{\rightarrow }\boldsymbol{\lambda }_{0}$
($N\rightarrow \infty $, $T$ fixed).%

\bigskip

\noindent \textbf{Theorem 4}. Suppose C1$^{\prime }$--C3$^{\prime }$, C4,
and C5 are satisfied and $\widetilde{\boldsymbol{H}}_{0}$ is negative
definite. Further assume the following condition is met:

\begin{description}
\item[C6$^{\prime }$:] the limit $\widetilde{\mathcal{I}}_{0}=\lim_{N%
\rightarrow \infty }N^{-1}\sum_{i}E\left[ \left( \partial \widetilde{l}%
_{i}\left( \boldsymbol{\lambda }_{0}\right) /\partial \boldsymbol{\lambda }%
\right) \left( \partial \widetilde{l}_{i}\left( \boldsymbol{\lambda }%
_{0}\right) /\partial \boldsymbol{\lambda }\right) ^{\prime }\right] $
exists and is positive definite.
\end{description}

\noindent Then $\sqrt{N}\left( \widehat{\boldsymbol{\lambda }}-\boldsymbol{%
\lambda }_{0}\right) \overset{d}{\rightarrow }\mathcal{N}\left( \mathbf{0},%
\widetilde{\boldsymbol{H}}_{0}^{-1}\widetilde{\mathcal{I}}_{0}\widetilde{%
\boldsymbol{H}}_{0}^{-1}\right) $ ($N\rightarrow \infty $, $T$ fixed).%

\bigskip

\noindent \textbf{Proof}. For proofs of Theorems 3 and 4, see Appendix D.%

The linear projection of $\Delta y_{i,-j+2}$ on 1 and $\boldsymbol{x}_{i}$
guarantees the residual in this linear projection is uncorrelated with the
elements of $\Delta \boldsymbol{X}_{i}$. This is a critical condition for
consistent differenced QML estimation. But this condition is also met if we
instead used the linear projection of $\Delta y_{i,-j+2}$ on 1 and $\Delta 
\boldsymbol{x}_{i}$, where $\Delta \boldsymbol{x}_{i}$ is a vector
consisting of the distinct elements of $\Delta \boldsymbol{X}_{i}$. The
latter approach generalizes an estimator studied by Hsiao et al. (2002).
Hsiao et al. (2002) studied differenced ML\ estimation of a dynamic panel
data model while assuming $p=1$, individual specific effects, and
uncorrelated and conditionally homoskedastic $v_{it}$s. Moreover, Hsiao et
al. (2002) also imposed restrictions on how the regressors are generated.
Furthermore, Hsiao et al. (2002) noted that the likelihood satisfies
standard regularity conditions, and therefore the ML estimator is consistent
and asymptotically normal. However, that conclusion follows from ML theory
assuming the log-likelihood is correctly specified. The analysis in this
section provides weaker conditions that imply the differenced ML estimator
proposed by Hsiao et al. (2002) is consistent and asymptotically normal (for 
$N\rightarrow \infty $, $T$ fixed). Specifically, the log-likelihood can be
misspecified and the $v_{it}$s can be conditionally heteroskedastic.  Moreover, all that is required of the
elements of $\boldsymbol{x}_{it}$ is that they be uncorrelated with the $%
v_{it}$s and that the linear projection of $\Delta y_{i1}$ on 1 and $\Delta 
\boldsymbol{x}_{i}$ does not depend on $i$.

\section{Computation}\label{computation}

If the error variance-covariance matrix is unrestricted, QML estimates can
be easily computed using iterated feasible generalized least squares. Consider, for example, calculating QML estimates of the elements of 
$\Omega _{0}$ and $\boldsymbol{\gamma }_{0}$. These estimates can be
calculated by iterating back and forth between fitting $\Omega _{0}$ and
fitting $\boldsymbol{\gamma }_{0}$. Specifically, $L_{N}\left( \cdot \right) 
$ is maximized with respect to the elements of $\Omega $, conditional on the
current fit of the regression parameters, say $\boldsymbol{\gamma }^{c}$, by
the fit $\Omega ^{+}=\sum_{i=1}^{N}\boldsymbol{u}_{i}\left( \boldsymbol{%
\gamma }^{c}\right) \boldsymbol{u}_{i}\left( \boldsymbol{\gamma }^{c}\right)
^{\prime }/N$. And, after $\Omega ^{+}$ is obtained, $L_{N}\left( \cdot
\right) $ is then maximized with respect to $\boldsymbol{\gamma }$,
conditional on $\Omega =\Omega ^{+}$, which gives the feasible generalized
least squares (FGLS) fit: 
\begin{equation}
\boldsymbol{\gamma }^{+}=\left( \sum_{i=1}^{N}\boldsymbol{W}_{i}^{\prime
}\left( \Omega ^{+}\right) ^{-1}\boldsymbol{W}_{i}\right) ^{-1}\sum_{i=1}^{N}%
\boldsymbol{W}_{i}^{\prime }\left( \Omega ^{+}\right) ^{-1}\boldsymbol{y}%
_{i}.  \label{fgls0}
\end{equation}%
This fit is then made the current fit, $\boldsymbol{\gamma }^{c}$, and new
fits $\Omega ^{+}$ and $\boldsymbol{\gamma }^{+}$ are calculated again, and
so on, until the sequence of fitted values converges. Calculating QML
estimates of $\boldsymbol{\lambda }_{0}$ and $\Upsilon _{0}$, based on
differenced observations, is similar when $\Upsilon _{0}$ is unrestricted.%

Although it is easy to calculate estimates by iterating back and forth
between fitting $\Omega _{0}$ and fitting $\boldsymbol{\gamma }_{0}$, or
between fitting $\boldsymbol{\lambda }_{0}$ and $\Upsilon _{0}$, this
approach implies that the number of free parameters being fitted in either $%
\Omega _{0}$ or $\Upsilon _{0}$ increases with $T$ at the rate $T^{2}$
increases. This fact, in turn, suggests that, if $T$ is not quite small, the
sampling performance of a QML estimator that does not impose valid
restrictions on $\Omega _{0}$ or $\Upsilon _{0}$ will be poor compared to
that of a QML estimator that does rely on valid restrictions.%

Unfortunately, maximizing the likelihood for differenced observations when
restrictions on $\Upsilon _{0}$ are imposed is tractable only for a highly
specialized case. Specifically, we must assume $p=1$, $e_{it}$ is given by
the error-components model in (\ref{EC}), the $v_{it}$s are uncorrelated and
unconditionally homoskedastic, and the regressors in $\boldsymbol{x}_{it}$
are strictly exogenous with respect to the $v_{it}$s. Further assume $\Delta
y_{i1}$ is generated by the same process generating $\Delta y_{it}$ for $%
t\geq 2$. Then it is easy to show that the error variance-covariance matrix
is $\Upsilon _{0}=\sigma _{0}^{2}\Phi _{0}$, 
\begin{equation}
\Phi _{0}=\left( 
\begin{array}{ccccc}
\phi _{0} & -1 & 0 & \cdots & 0 \\ 
-1 & 2 & -1 & \cdots & 0 \\ 
0 & -1 & 2 & \ddots & \vdots \\ 
\vdots & \vdots & \ddots & \ddots & -1 \\ 
0 & 0 & \cdots & -1 & 2%
\end{array}%
\right)  \label{Phi_0}
\end{equation}%
(cf Hsiao et al. 2002, p. 110, Eq. (3.2)). Moreover, the determinant $%
\left\vert \sigma _{0}^{2}\Phi _{0}\right\vert $ equals $\sigma _{0}^{2T}%
\left[ 1+T\left( \phi _{0}-1\right) \right] $ (see, e.g., Hsiao et al 2002,
p. 111, Eq. (3.7)). From this determinant we see that, in order to ensure a
positive definite fitted value for $\sigma _{0}^{2}\Phi _{0}$, we must
search over values of $\phi $ satisfying $\phi >1-1/T$. This restriction is
guaranteed if we set $\varpi =\ln \left( \phi -1+1/T\right) $ and maximize
the log-likelihood%
\begin{equation*}
const-\frac{NT}{2}\ln \left( \sigma ^{2}\right) -\frac{N\varpi }{2}-\frac{1}{%
2\sigma ^{2}}\sum_{i=1}^{N}\widetilde{\boldsymbol{u}}_{i}\left( \boldsymbol{%
\eta }\right) ^{\prime }\Phi ^{-1}\widetilde{\boldsymbol{u}}_{i}\left( 
\boldsymbol{\eta }\right)
\end{equation*}%
with respect to $\boldsymbol{\eta }$, $\sigma ^{2}$, and $\varpi $. Here $%
\Phi $ has $\exp \left( \varpi \right) +1-1/T$ in its first row, first
column and everywhere else is the same as $\Phi _{0}$ in (\ref{Phi_0}).%

Maximizing the log-likelihood for differenced QML estimation becomes much
more complicated if the $v_{it}$s are time-series heteroskedastic or $p>1$.
On the other hand, the ease with which levels QML estimates can be
calculated is not affected by the size of $p$ nor by whether or not the $%
v_{it}$s are time-series heteroskedastic. The remainder of this section is
devoted to describing an ECME algorithm that can be applied to calculate
levels QML estimates for arbitrary $p$ and for an error variance-covariance
matrix given by $\Omega _{0}=\sigma _{a0}^{2}\boldsymbol{\iota \iota }%
^{\prime }+\Sigma _{0}$, with $\Sigma _{0}=$ diag$\left( \sigma
_{01}^{2},\ldots ,\sigma _{0T}^{2}\right) $.%

The ECME algorithm relies on conditional or constrained maximization (CM) of
either an imputed log-likelihood, based on augmented data, or the
log-likelihood based on the observed data. In the present application, the
observed data are $\boldsymbol{y}=\left( \boldsymbol{y}_{1}^{\prime },\ldots
,\boldsymbol{y}_{N}^{\prime }\right) ^{\prime }$, while the augmented data
consists of $\boldsymbol{y}$ and $\boldsymbol{a}=\left( a_{1},\ldots
,a_{N}\right) ^{\prime }$.\footnote{%
For the purposes of deriving the imputed log-likelihood and the actual
log-likelihood, the variables in $\boldsymbol{z}=\left( \boldsymbol{z}%
_{1}^{\prime },\ldots ,\boldsymbol{z}_{N}^{\prime }\right) ^{\prime }$ are
treated as fixed.} The imputed log-likelihood is built during the
expectation (E) step by taking the conditional expectation of the
log-likelihood for the augmented data given the observed data, while
treating the current fit of the parameters $\boldsymbol{\psi }^{c}$ as the
parameters of the conditional distribution.\footnote{%
Liu and Rubin (1994) describe the properties of the ECME algorithm. For
applications of it to panel data see Phillips (2004, 2012).}%

Applying the ECME algorithm to an error-components model for which $\Omega
_{0}=\sigma _{a0}^{2}\boldsymbol{\iota \iota }^{\prime }+\Sigma _{0}$, with $%
\Sigma _{0}=$ diag$\left( \sigma _{01}^{2},\ldots ,\sigma _{0T}^{2}\right) $%
, leads to the following E and CM steps:

\textit{E-step}: Let $\left( \sigma _{a}^{2}\right) ^{c}$, $\boldsymbol{%
\gamma }^{c}$, and $\Omega ^{c}=\left( \sigma _{a}^{2}\right) ^{c}%
\boldsymbol{\iota \iota }^{\prime }+\Sigma ^{c}$, with $\Sigma ^{c}=$ diag$%
\left( \left( \sigma _{1}^{2}\right) ^{c},\ldots ,\left( \sigma
_{T}^{2}\right) ^{c}\right) $, denote the current fits of $\sigma _{a0}^{2}$%
, $\boldsymbol{\gamma }_{0}$, and $\Omega _{0}$. Compute the conditional
mean and variance of $a_{i}$ given $\boldsymbol{y}_{i}$ evaluated at the
current fit of the parameters. These are $a_{i}^{c}=\left( \sigma
_{a}^{2}\right) ^{c}\boldsymbol{\iota \,}^{\prime }\left( \Omega ^{c}\right)
^{-1}\boldsymbol{u}_{i}\left( \boldsymbol{\gamma }^{c}\right) $ and $%
\upsilon _{a}^{c}=\left( \sigma _{a}^{2}\right) ^{c}\left[ 1-\left( \sigma
_{a}^{2}\right) ^{c}\boldsymbol{\iota }^{\prime }\left( \Omega ^{c}\right)
^{-1}\boldsymbol{\iota }\right] $, respectively (see, e.g., Greene 2012,
Theorem B.7, pp. 1041-1042). Then the imputed log-likelihood is%
\begin{eqnarray*}
Q\left( \boldsymbol{\psi \,};\boldsymbol{\psi }^{c}\right) &=&const-\frac{N}{%
2}\left( \ln \sigma _{a}^{2}+\sum_{t=1}^{T}\ln \sigma _{t}^{2}\right) -\frac{%
1}{2\sigma _{a}^{2}}\sum_{i=1}^{N}\left( a_{i}^{c}\right) ^{2}-\frac{N}{%
2\sigma _{a}^{2}}\upsilon _{a}^{c} \\
&&-\frac{1}{2}\sum_{i=1}^{N}\left[ \boldsymbol{u}_{i}\left( \boldsymbol{%
\gamma }\right) -\boldsymbol{\iota \,}a_{i}^{c}\right] ^{\prime }\Sigma ^{-1}%
\left[ \boldsymbol{u}_{i}\left( \boldsymbol{\gamma }\right) -\boldsymbol{%
\iota \,}a_{i}^{c}\right] -\frac{N}{2}\boldsymbol{\iota }\,^{\prime }\Sigma
^{-1}\boldsymbol{\iota \,}\upsilon _{a}^{c}.
\end{eqnarray*}%

\textit{CM-step} 1: Maximize $Q\left( \boldsymbol{\cdot };\boldsymbol{%
\psi }^{c}\right) $ with respect to $\boldsymbol{\omega }=\left( \sigma
_{a}^{2},\sigma _{1}^{2},\ldots ,\sigma _{T}^{2}\right) ^{\prime }$ subject
to the constraint $\boldsymbol{\gamma }=\boldsymbol{\gamma }^{c}$. This step
yields $\left( \sigma _{a}^{2}\right) ^{+}=\upsilon
_{a}^{c}+\sum_{i=1}^{N}\left( a_{i}^{c}\right) ^{2}/N$ and%
\begin{equation}
\left( \sigma _{t}^{2}\right) ^{+}=\upsilon _{a}^{c}+\frac{1}{N}%
\sum_{i=1}^{N}\left[ u_{it}\left( \boldsymbol{\gamma }^{c}\right) -a_{i}^{c}%
\right] ^{2}\text{ \ \ \ \ \ \ \ \ }t=1,\ldots ,T.  \label{vart}
\end{equation}%

\textit{CM-step} 2: Maximize the actual log-likelihood $%
\sum_{i=1}^{N}l_{i}\left( \cdot \right) $ with respect to $\boldsymbol{%
\gamma }$ subject to the constraint $\boldsymbol{\omega }=\boldsymbol{\omega 
}^{+}$, where $\boldsymbol{\omega }^{+}=\left( \left( \sigma _{a}^{2}\right)
^{+},\left( \sigma _{1}^{2}\right) ^{+},\ldots ,\left( \sigma
_{T}^{2}\right) ^{+}\right) ^{\prime }$. This step gives the FGLS fit in Eq.
(\ref{fgls0}) with $\Omega ^{+}=\left( \sigma _{a}^{2}\right) ^{+}%
\boldsymbol{\iota \iota }^{\prime }+\Sigma ^{+}$ and $\Sigma ^{+}=$ diag$%
\left( \left( \sigma _{1}^{2}\right) ^{+},\ldots ,\left( \sigma
_{T}^{2}\right) ^{+}\right) $.

After the new fits of the parameters are obtained, they become the current
fits, and the preceding steps are repeated, until convergence.%

Unlike some other algorithms, the ECME fitted values for the error variance
components are guaranteed to be non-negative. But this advantage can lead to
another complication. Specifically, EM-like algorithms --- including the
ECME algorithm --- can be excruciatingly slow to converge, and, when
calculating estimates of error-components models, the rate of convergence
can slow when the sequence of the fitted variance of the individual-specific
effect gets close to zero (see Meng and van Dyk 1998). Moreover, there is
always the possibility that the error-components model in (\ref{EC}) is
inappropriate; specifically, there may be no individual-specific effects. In
this case, we have $\sigma _{c0}^{2}=0$, where $\sigma _{c0}^{2}=var\left(
c_{i}\right) $, and $\sigma _{a0}^{2}=0$, and consequently the sequence of
fitted values for $\sigma _{a0}^{2}$ can approach zero. Furthermore, even if 
$\sigma _{c0}^{2}$ is positive and large, $\sigma _{a0}^{2}$ can be small,
for the control function $\mu _{0}+\boldsymbol{z}_{i}^{\prime }\boldsymbol{%
\theta }_{0}$ is the best linear predictor of $c_{i}$ based on $\boldsymbol{z%
}_{i}$, and if that predictor is accurate, then $\sigma _{a0}^{2}$ can be
near zero. If so, the sequence of fitted values for $\sigma _{a0}^{2}$ can
get close to zero.

As a practical matter, however, given $\Omega _{0}=\sigma _{a0}^{2}%
\boldsymbol{\iota \iota }^{\prime }+\Sigma _{0}$, with $\Sigma _{0}=$ diag$%
\left( \sigma _{01}^{2},\ldots ,\sigma _{0T}^{2}\right) $, then, when the
fitted value for $\sigma _{a0}^{2}$ is near zero, the fitted value $%
\boldsymbol{\gamma }^{+}$ in (\ref{fgls0}) differs little from the weighted
least squares fit $\left( \sum_{i=1}^{N}\boldsymbol{W}_{i}^{\prime }\left(
\Sigma ^{+}\right) ^{-1}\boldsymbol{W}_{i}\right) ^{-1}\sum_{i=1}^{N}%
\boldsymbol{W}_{i}^{\prime }\left( \Sigma ^{+}\right) ^{-1}\boldsymbol{y}_{i}
$, which is obtained by setting $\left( \sigma _{a}^{2}\right) ^{+}=0$.
Furthermore, once $\left( \sigma _{a}^{2}\right) ^{+}$ is set to zero, all
subsequent fitted values for $\sigma _{a0}^{2}$ will be zero. Also, when $%
\left( \sigma _{a}^{2}\right) ^{c}=0$, Eq. (\ref{vart}) simplifies to $%
\left( \sigma _{t}^{2}\right) ^{+}=\sum_{i=1}^{N}u_{it}\left( \boldsymbol{%
\gamma }^{c}\right) ^{2}/N$. Thus, if $\left( \sigma _{a}^{2}\right) ^{+}$
is set to zero, convergence is rapid. Consequently, the ECME algorithm for
computing level QML estimates will generally converge at a robust rate if,
as part of the convergence criterion, the size of the fitted value for $%
\sigma _{a0}^{2}$ is evaluated and $\left( \sigma _{a}^{2}\right) ^{+}$ is
set to zero should it become sufficiently small.\footnote{%
For example, the fitted value of $\sigma _{a}^{2}$ might be set to zero when
the fitted value for the average correlation coefficient, say $\overline{%
\rho }$, is small, where $\overline{\rho }=2\sum_{s=1}^{T-1}\sum_{t>s}^{T}%
\rho _{st}/\left[ T\left( T-1\right) \right] $, with $\rho _{st}=\sigma
_{a0}^{2}/\left[ \left( \sigma _{a0}^{2}+\sigma _{0s}^{2}\right) \left(
\sigma _{a0}^{2}+\sigma _{0t}^{2}\right) \right] ^{1/2}$. This criterion was
used to obtain the results for the levels QML estimator provided in Section %
\ref{results}. In particular, the fitted value of $\sigma _{a}^{2}$ was set
to zero when the fitted value of $\rho $ fell below 0.01.}

\section{Monte Carlo Experiments}

\subsection{Design}

In order to assess the finite sampling properties of QML estimators
described in Section \ref{computation}, Monte Carlo experiments were
conducted. For all of the experiments, observations on the dependent
variable $y_{it}$ were generated according to the model%
\begin{equation*}
y_{it}=\delta _{0}y_{i,t-1}+0.5x_{it}+c_{i}+v_{it}\text{ \ \ \ \ \ \ \ \ \ \
\ \ \ }\left( t=-t_{0}+1,\ldots ,T,\text{ \ }i=1,\ldots ,N\right) ,
\end{equation*}%
with $y_{i,-t_{0}}=0$. The values for $\delta _{0}$ considered were 0, 0.2,
0.4, 0.6, 0.8, and 0.9. Moreover, the $x_{it}$s were generated according to
the autoregressive process%
\begin{equation*}
x_{it}=0.5+0.5x_{i,t-1}+\xi _{it}\text{ \ \ \ \ \ \ \ \ \ \ \ \ \ }\left(
t=-t_{0}+1,\ldots ,T,\text{ \ }i=1,\ldots ,N\right) .
\end{equation*}%
The starting value $x_{i,-t_{0}}$ was set equal to $5+10\xi _{i,-t_{0}}$ and
the $\xi _{it}$s were generated as independent uniform random variates with
mean zero and variance one. Furthermore, two values for $t_{0}$ were
considered: $t_{0}=1$ and $t_{0}=50$. For $t_{0}=50$, the time series for $%
x_{it}$ and $y_{it}$ were essentially stationary, whereas for $t_{0}=1$ they
were nonstationary.

As for the $v_{it}$s, they were generated as $v_{it}=x_{it}\left( \epsilon
_{it}-5\right) /\sqrt{10}$, with $\epsilon _{it}$ a chi-square random
variate with five degress of freedom. The variate $\left( \epsilon
_{it}-5\right) /\sqrt{10}$ has an asymmetric distribution about zero with a
variance of one. Moreover, because the $\epsilon _{it}$s were generated
independently of one another and of the $x_{it}$s, the $v_{it}$s were
uncorrelated but conditionally heteroskedastic. However, the $v_{it}$s were unconditionally homoskedastic for $t \geq 1$when $t_{0}$ was
set to 50, for in this case the $x_{it}$s were essentially stationary by the time $t = 1$. On the other hand, for $t_{0}=1$%
, the $x_{it}$s had insufficient time to become approximately stationary by the time $t = 1$. Hence, in this case, the $v_{it}$s were not only conditionally heteroskedastic, they were also
unconditionally time-series heteroskedastic for $t \geq 1$.%

\sloppy The heterogeneity component, $c_{i}$, was generated as $%
c_{i}=\sum_{t=0}^{T}\ln |x_{it}|/\left( T+1\right) +\sigma _{\zeta }\left(
\zeta _{i}-5\right) /\sqrt{10}$, with $\zeta _{i}$ a chi-square random
variate with five degress of freedom.  Furthermore, the parameter $\sigma _{\zeta }$ was set to
either one or four. This specification for $c_{i}$ induced correlation
between $c_{i}$ and the $x_{it}$s. Moreover, both $c_{i}$ and $v_{it}$,
conditional on the $x_{it}$s, had non-normal asymmetric distributions,
implying that, conditional on the $x_{it}$s, the error $e_{it}=c_{i}+v_{it}$
came from a non-normal asymmetric distribution.%

\fussy After a sample was generated, the start up observations were discarded so
that QML estimation was based on $\left( x_{i1},y_{i1}\right) ,\ldots
,\left( x_{iT},y_{iT}\right) $ and $y_{i0}$ ($i=1,\ldots ,N$), while GMM
estimation was based on $\left( x_{i0},y_{i0}\right) ,\ldots ,\left(
x_{iT},y_{iT}\right) $. Furthermore, $T$ was set to ten, and $N$ was set to
200. Finally, for each combination of parameters, 5,000 independent samples were
generated.

\subsection{Estimators\label{estimators}}

The finite sample properties of levels and differenced QML estimators
were compared to each other and to two well-known GMM estimators. The GMM
estimators considered were the differenced GMM estimator proposed by
Arellano and Bond (1991) (denoted DGMM) and the system GMM estimator
suggested by Blundell and Bond (1998) (SGMM).%

Three QML estimators were considered. Results
are provided for levels QML (LQML) while relying on the
structured variance-covariance matrix $\Omega _{0}=\sigma _{a0}^{2}%
\boldsymbol{\iota \iota }^{\prime }+\Sigma _{0}$ with $\Sigma _{0}=$ diag$%
\left( \sigma _{01}^{2},\ldots ,\sigma _{0T}^{2}\right) $. For this case,
estimates were calculated with the ECME algorithm. Differenced QML estimates
were also calculated. As noted in Section \ref{computation}, computing
differenced QML estimates via gradient methods is complicated if we model
the $v_{it}$s as time-series heteroskedastic. For this reason, results are
only provided for differenced QML estimates that restrict the $v_{it}$s to
be uncorrelated and unconditionally homoskedastic. Because we can use either
a linear projection of $\Delta y_{i1}$ on 1 and $\Delta \boldsymbol{x}_{i}$
or a linear projection of $\Delta y_{i1}$ on 1 and $\boldsymbol{x}_{i}$,
results for both choices are reported and are denoted by DQML$_{\Delta 
\boldsymbol{x}}$ and DQML$_{\boldsymbol{x}}$.%

\subsection{Results\label{results}}

\subsubsection{Stationary Designs}

This section provides results for designs for which the generated variables
were approximately stationary ($t_{0}=50)$. Table 1 provides estimates of
finite sample bias and root mean squared error for the panel data GMM and  QML estimators for stationary designs with $\sigma _{\zeta }=1$ and $\sigma _{\zeta }=4$.

The evidence in Table 1 shows that the QML estimators --- LQML, DQML$_{%
\boldsymbol{x}}$, and DQML$_{\Delta \boldsymbol{x}}$ --- generally have
neglible finite sample bias, and, consequently, their root mean squared
errors are significantly smaller than that of the GMM estimators, which have
non-neglible finite sample bias. Moreover, for most designs, whether one
uses DQML$_{\boldsymbol{x}}$ or DQML$_{\Delta \boldsymbol{x}}$ does not
matter much; they have similar finite sample bias and root mean squared
error. The exception is when $\delta _{0}=0.9$. For highly persistent designs,
DQML$_{\boldsymbol{x}}$ outperforms DQML$_{\Delta \boldsymbol{x}}$. But
among the QML estimators, the levels QML estimator (LQML) is --- in terms of
root mean squared error --- best.

The system GMM estimator was introduced as a response to the poor sampling
performance of the differenced GMM estimator when $\delta _{0} $ is near one.
Blundell and Bond (1998) showed that the system GMM estimator will perform
better than the differenced GMM estimator in this case, and it does indeed
have smaller bias and root mean squared error than the differenced GMM
estimator for $\delta _{0} $ near one and $\sigma _{\zeta }=1$. However, surprisingly, its sampling
performance is worse --- often much worse --- than that of the differenced
GMM estimator for $\delta _{0} $ not near one. Furthermore, when $\sigma _{\zeta }=4$, the system GMM estimator has substantial bias even when $\delta _{0} $ is near one. Bun and Windmeijer (2010) provide an explanation for this result. They note the system GMM estimator may
suffer from a weak instrument problem when the variance of the
individual-specific effect is large relative to the variance of the
idiosyncratic error. The sampling performance of the QML estimators, on the
other hand, are unaffected by the relative size of the individual-specific
effect variance versus the idiosyncratic error variance.

 \vspace{.75cm}
\begin{center}
\textbf{Table 1:} Finite sample characteristics of estimators of $\delta
_{0}$ for $t=50$.
\end{center}
\begin{center}
\begin{tabular}{llcccccc}
\hline
&  &  &  &  &  &  &  \\ 
&  & \multicolumn{6}{c}{$\delta _{0}$} \\ 
&  & $0.0$ & $0.2$ & $0.4$ & $0.6$ & $0.8$ & $0.9$ \\ \hline
& \multicolumn{1}{c}{} &  &  & \multicolumn{1}{r}{} &  &  &  \\ 
$\sigma _{\zeta }=1$ &  &  &  &  &  &  &  \\ \cline{1-1}
&  &  &  &  &  &  &  \\ 
DGMM & bias & \multicolumn{1}{r}{$-0.0112$} & \multicolumn{1}{r}{$-0.0147$}
& \multicolumn{1}{r}{$-0.0223$} & \multicolumn{1}{r}{$-0.0318$} & 
\multicolumn{1}{r}{$-0.0533$} & \multicolumn{1}{r}{$-0.0784$} \\ 
& \textit{rmse} & \multicolumn{1}{r}{$\mathit{0.0322}$} & \multicolumn{1}{r}{%
$\mathit{0.0345}$} & \multicolumn{1}{r}{$\mathit{0.0391}$} & 
\multicolumn{1}{r}{$\mathit{0.0455}$} & \multicolumn{1}{r}{$\mathit{0.0637}$}
& \multicolumn{1}{r}{$\mathit{0.0875}$} \\ 
&  &  &  &  &  &  & \multicolumn{1}{r}{} \\ 
SGMM & bias & \multicolumn{1}{r}{$-0.0419$} & \multicolumn{1}{r}{$-0.0547$}
& \multicolumn{1}{r}{$-0.0679$} & \multicolumn{1}{r}{$-0.0741$} & 
\multicolumn{1}{r}{$-0.0498$} & \multicolumn{1}{r}{$-0.0106$} \\ 
& \textit{rmse} & \multicolumn{1}{r}{$\mathit{0.0539}$} & \multicolumn{1}{r}{%
$\mathit{0.0662}$} & \multicolumn{1}{r}{$\mathit{0.0785}$} & 
\multicolumn{1}{r}{$\mathit{0.0850}$} & \multicolumn{1}{r}{$\mathit{0.0641}$}
& \multicolumn{1}{r}{$\mathit{0.0332}$} \\ 
&  &  &  &  &  &  & \multicolumn{1}{r}{} \\ 
LQML & bias & \multicolumn{1}{r}{$0.0003$} & \multicolumn{1}{r}{$0.0005$} & 
\multicolumn{1}{r}{$-0.0005$} & \multicolumn{1}{r}{$-0.0004$} & 
\multicolumn{1}{r}{$-0.0001$} & \multicolumn{1}{r}{$-0.0049$} \\ 
& \textit{rmse} & \multicolumn{1}{r}{$\mathit{0.0275}$} & \multicolumn{1}{r}{%
$\mathit{0.0282}$} & \multicolumn{1}{r}{$\mathit{0.0273}$} & 
\multicolumn{1}{r}{$\mathit{0.0265}$} & \multicolumn{1}{r}{$\mathit{0.0271}$}
& \multicolumn{1}{r}{$\mathit{0.0269}$} \\ 
&  &  &  &  &  &  & \multicolumn{1}{r}{} \\ 
DQML$_{\boldsymbol{x}}$ & bias & \multicolumn{1}{r}{$0.0002$} & 
\multicolumn{1}{r}{$0.0005$} & \multicolumn{1}{r}{$-0.0005$} & 
\multicolumn{1}{r}{$-0.0002$} & \multicolumn{1}{r}{$0.0011$} & 
\multicolumn{1}{r}{$-0.0001$} \\ 
& \textit{rmse} & \multicolumn{1}{r}{$\mathit{0.0281}$} & \multicolumn{1}{r}{%
$\mathit{0.0287}$} & \multicolumn{1}{r}{$\mathit{0.0281}$} & 
\multicolumn{1}{r}{$\mathit{0.0276}$} & \multicolumn{1}{r}{$\mathit{0.0309}$}
& \multicolumn{1}{r}{$\mathit{0.0369}$} \\ 
& \multicolumn{1}{c}{} & \multicolumn{1}{r}{} & \multicolumn{1}{r}{} & 
\multicolumn{1}{r}{} &  & \multicolumn{1}{r}{} & \multicolumn{1}{r}{} \\ 
DQML$_{\Delta \boldsymbol{x}}$ & bias & \multicolumn{1}{r}{$0.0002$} & 
\multicolumn{1}{r}{$0.0005$} & \multicolumn{1}{r}{$-0.0005$} & 
\multicolumn{1}{r}{$-0.0002$} & \multicolumn{1}{r}{$0.0012$} & 
\multicolumn{1}{r}{$0.0008$} \\ 
& \textit{rmse} & \multicolumn{1}{r}{$\mathit{0.0281}$} & \multicolumn{1}{r}{%
$\mathit{0.0287}$} & \multicolumn{1}{r}{$\mathit{0.0280}$} & 
\multicolumn{1}{r}{$\mathit{0.0276}$} & \multicolumn{1}{r}{$\mathit{0.0310}$}
& \multicolumn{1}{r}{$\mathit{0.0389}$} \\ 

& \multicolumn{1}{c}{} &  &  & \multicolumn{1}{r}{} &  &  &  \\ 
$\sigma _{\zeta }=4$ &  &  &  &  &  &  &  \\ \cline{1-1}
&  &  &  &  &  &  &  \\ 
DGMM & bias & \multicolumn{1}{r}{$-0.0139$} & \multicolumn{1}{r}{$-0.0183$}
& \multicolumn{1}{r}{$-0.0248$} & \multicolumn{1}{r}{$-0.0366$} & 
\multicolumn{1}{r}{$-0.0616$} & \multicolumn{1}{r}{$-0.0809$} \\ 
& \textit{rmse} & \multicolumn{1}{r}{$\mathit{0.0346}$} & \multicolumn{1}{r}{%
$\mathit{0.0383}$} & \multicolumn{1}{r}{$\mathit{0.0425}$} & 
\multicolumn{1}{r}{$\mathit{0.0508}$} & \multicolumn{1}{r}{$\mathit{0.0721}$}
& \multicolumn{1}{r}{$\mathit{0.0899}$} \\ 
&  &  &  &  &  &  & \multicolumn{1}{r}{} \\ 
SGMM & bias & \multicolumn{1}{r}{$-0.0058$} & \multicolumn{1}{r}{$0.0057$} & 
\multicolumn{1}{r}{$0.0313$} & \multicolumn{1}{r}{$0.0736$} & 
\multicolumn{1}{r}{$0.1032$} & \multicolumn{1}{r}{$0.0769$} \\ 
& \textit{rmse} & \multicolumn{1}{r}{$\mathit{0.0416}$} & \multicolumn{1}{r}{%
$\mathit{0.0482}$} & \multicolumn{1}{r}{$\mathit{0.0630}$} & 
\multicolumn{1}{r}{$\mathit{0.0896}$} & \multicolumn{1}{r}{$\mathit{0.1086}$}
& \multicolumn{1}{r}{$\mathit{0.0780}$} \\ 
&  &  &  &  &  &  & \multicolumn{1}{r}{} \\ 
LQML & bias & \multicolumn{1}{r}{$-0.0006$} & \multicolumn{1}{r}{$-0.0001$}
& \multicolumn{1}{r}{$-0.0002$} & \multicolumn{1}{r}{$0.0003$} & 
\multicolumn{1}{r}{$-0.0011$} & \multicolumn{1}{r}{$-0.0038$} \\ 
& \textit{rmse} & \multicolumn{1}{r}{$\mathit{0.0278}$} & \multicolumn{1}{r}{%
$\mathit{0.0283}$} & \multicolumn{1}{r}{$\mathit{0.0278}$} & 
\multicolumn{1}{r}{$\mathit{0.0271}$} & \multicolumn{1}{r}{$\mathit{0.0277}$}
& \multicolumn{1}{r}{$\mathit{0.0276}$} \\ 
&  &  &  &  &  &  & \multicolumn{1}{r}{} \\ 
DQML$_{\boldsymbol{x}}$ & bias & \multicolumn{1}{r}{$-0.0007$} & 
\multicolumn{1}{r}{$-0.0001$} & \multicolumn{1}{r}{$-0.0001$} & 
\multicolumn{1}{r}{$0.0006$} & \multicolumn{1}{r}{$-0.0006$} & 
\multicolumn{1}{r}{$0.0003$} \\ 
& \textit{rmse} & \multicolumn{1}{r}{$\mathit{0.0281}$} & \multicolumn{1}{r}{%
$\mathit{0.0288}$} & \multicolumn{1}{r}{$\mathit{0.0285}$} & 
\multicolumn{1}{r}{$\mathit{0.0279}$} & \multicolumn{1}{r}{$\mathit{0.0303}$}
& \multicolumn{1}{r}{$\mathit{0.0366}$} \\ 
& \multicolumn{1}{c}{} & \multicolumn{1}{r}{} & \multicolumn{1}{r}{} & 
\multicolumn{1}{r}{} &  & \multicolumn{1}{r}{} & \multicolumn{1}{r}{} \\ 
DQML$_{\Delta \boldsymbol{x}}$ & bias & \multicolumn{1}{r}{$-0.0007$} & 
\multicolumn{1}{r}{$-0.0001$} & \multicolumn{1}{r}{$-0.0001$} & 
\multicolumn{1}{r}{$0.0006$} & \multicolumn{1}{r}{$-0.0006$} & 
\multicolumn{1}{r}{$0.0010$} \\ 
& \textit{rmse} & \multicolumn{1}{r}{$\mathit{0.0281}$} & \multicolumn{1}{r}{%
$\mathit{0.0288}$} & \multicolumn{1}{r}{$\mathit{0.0284}$} & 
\multicolumn{1}{r}{$\mathit{0.0278}$} & \multicolumn{1}{r}{$\mathit{0.0302}$}
& \multicolumn{1}{r}{$\mathit{0.0386}$} \\ 
&  &  &  &  &  &  & \multicolumn{1}{r}{} \\ \hline
\end{tabular}%
\bigskip 
\end{center}

\subsubsection{Nonstationary Designs\label{nonstat_designs}}

Table 2 provides finite sample bias and root mean squared error
estimates for nonstationary designs. For these designs $t_{0}=1$, and,
therefore, for each cross section, the time series began in the immediate
past.

\vspace{.4cm}
\begin{center}
\noindent \textbf{Table 2:} Finite sample characteristics of estimators of $\delta
_{0}$ for $t=1$.
\end{center}

\begin{center}
\begin{tabular}{llcccccc}
\hline
&  &  &  &  &  &  &  \\ 
&  & \multicolumn{6}{c}{$\delta _{0}$} \\
&  & $0.0$ & $0.2$ & $0.4$ & $0.6$ & $0.8$ & $0.9$ \\ \hline
& \multicolumn{1}{c}{} &  &  & \multicolumn{1}{r}{} &  &  &  \\ 
$\sigma _{\zeta }=1$ &  &  &  &  &  &  &  \\ \cline{1-1}
&  &  &  &  &  &  &  \\ 
DGMM & bias & \multicolumn{1}{r}{$-0.0042$} & \multicolumn{1}{r}{$-0.0050$}
& \multicolumn{1}{r}{$-0.0086$} & \multicolumn{1}{r}{$-0.0129$} & 
\multicolumn{1}{r}{$-0.0257$} & \multicolumn{1}{r}{$-0.0349$} \\ 
& \textit{rmse} & \multicolumn{1}{r}{$\mathit{0.0329}$} & \multicolumn{1}{r}{%
$\mathit{0.0340}$} & \multicolumn{1}{r}{$\mathit{0.0333}$} & 
\multicolumn{1}{r}{$\mathit{0.0332}$} & \multicolumn{1}{r}{$\mathit{0.0408}$}
& \multicolumn{1}{r}{$\mathit{0.0463}$} \\ 
&  &  &  &  &  &  & \multicolumn{1}{r}{} \\ 
SGMM & bias & \multicolumn{1}{r}{$-0.0179$} & \multicolumn{1}{r}{$-0.0226$}
& \multicolumn{1}{r}{$-0.0262$} & \multicolumn{1}{r}{$-0.0208$} & 
\multicolumn{1}{r}{$0.0097$} & \multicolumn{1}{r}{$0.0581$} \\ 
& \textit{rmse} & \multicolumn{1}{r}{$\mathit{0.0382}$} & \multicolumn{1}{r}{%
$\mathit{0.0422}$} & \multicolumn{1}{r}{$\mathit{0.0439}$} & 
\multicolumn{1}{r}{$\mathit{0.0390}$} & \multicolumn{1}{r}{$\mathit{0.0317}$}
& \multicolumn{1}{r}{$\mathit{0.0630}$} \\ 
&  &  &  &  &  &  & \multicolumn{1}{r}{} \\ 
LQML & bias & \multicolumn{1}{r}{$-0.0005$} & \multicolumn{1}{r}{$0.0001$} & 
\multicolumn{1}{r}{$-0.0012$} & \multicolumn{1}{r}{$-0.0010$} & 
\multicolumn{1}{r}{$-0.0010$} & \multicolumn{1}{r}{$0.0001$} \\ 
& \textit{rmse} & \multicolumn{1}{r}{$\mathit{0.0248}$} & \multicolumn{1}{r}{%
$\mathit{0.0250}$} & \multicolumn{1}{r}{$\mathit{0.0233}$} & 
\multicolumn{1}{r}{$\mathit{0.0218}$} & \multicolumn{1}{r}{$\mathit{0.0221}$}
& \multicolumn{1}{r}{$\mathit{0.0226}$} \\ 
&  &  &  &  &  &  & \multicolumn{1}{r}{} \\ 
DQML$_{\boldsymbol{x}}$ & bias & \multicolumn{1}{r}{$-0.0021$} & 
\multicolumn{1}{r}{$-0.0029$} & \multicolumn{1}{r}{$-0.0068$} & 
\multicolumn{1}{r}{$-0.0127$} & \multicolumn{1}{r}{$-0.0291$} & 
\multicolumn{1}{r}{$-0.0400$} \\ 
& \textit{rmse} & \multicolumn{1}{r}{$\mathit{0.0328}$} & \multicolumn{1}{r}{%
$\mathit{0.0332}$} & \multicolumn{1}{r}{$\mathit{0.0320}$} & 
\multicolumn{1}{r}{$\mathit{0.0315}$} & \multicolumn{1}{r}{$\mathit{0.0409}$}
& \multicolumn{1}{r}{$\mathit{0.0493}$} \\ 
& \multicolumn{1}{c}{} & \multicolumn{1}{r}{} & \multicolumn{1}{r}{} & 
\multicolumn{1}{r}{} &  & \multicolumn{1}{r}{} & \multicolumn{1}{r}{} \\ 
DQML$_{\Delta \boldsymbol{x}}$ & bias & \multicolumn{1}{r}{$-0.0021$} & 
\multicolumn{1}{r}{$-0.0029$} & \multicolumn{1}{r}{$-0.0068$} & 
\multicolumn{1}{r}{$-0.0128$} & \multicolumn{1}{r}{$-0.0295$} & 
\multicolumn{1}{r}{$-0.0417$} \\ 
& \textit{rmse} & \multicolumn{1}{r}{$\mathit{0.0328}$} & \multicolumn{1}{r}{%
$\mathit{0.0332}$} & \multicolumn{1}{r}{$\mathit{0.0320}$} & 
\multicolumn{1}{r}{$\mathit{0.0315}$} & \multicolumn{1}{r}{$\mathit{0.0411}$}
& \multicolumn{1}{r}{$\mathit{0.0504}$} \\ 

& \multicolumn{1}{c}{} &  &  & \multicolumn{1}{r}{} &  &  &  \\ 
$\sigma _{\zeta }=4$ &  &  &  &  &  &  &  \\ \cline{1-1}
&  &  &  &  &  &  &  \\ 
DGMM & bias & \multicolumn{1}{r}{$-0.0070$} & \multicolumn{1}{r}{$-0.0118$}
& \multicolumn{1}{r}{$-0.0193$} & \multicolumn{1}{r}{$-0.0317$} & 
\multicolumn{1}{r}{$-0.0110$} & \multicolumn{1}{r}{$-0.0043$} \\ 
& \textit{rmse} & \multicolumn{1}{r}{$\mathit{0.0357}$} & \multicolumn{1}{r}{%
$\mathit{0.0390}$} & \multicolumn{1}{r}{$\mathit{0.0455}$} & 
\multicolumn{1}{r}{$\mathit{0.0480}$} & \multicolumn{1}{r}{$\mathit{0.0193}$}
& \multicolumn{1}{r}{$\mathit{0.0100}$} \\ 
&  &  &  &  &  &  & \multicolumn{1}{r}{} \\ 
SGMM & bias & \multicolumn{1}{r}{$-0.0055$} & \multicolumn{1}{r}{$-0.0163$}
& \multicolumn{1}{r}{$0.0942$} & \multicolumn{1}{r}{$0.2512$} & 
\multicolumn{1}{r}{$0.2636$} & \multicolumn{1}{r}{$0.1914$} \\ 
& \textit{rmse} & \multicolumn{1}{r}{$\mathit{0.0385}$} & \multicolumn{1}{r}{%
$\mathit{0.0489}$} & \multicolumn{1}{r}{$\mathit{0.1110}$} & 
\multicolumn{1}{r}{$\mathit{0.2568}$} & \multicolumn{1}{r}{$\mathit{0.2638}$}
& \multicolumn{1}{r}{$\mathit{0.1915}$} \\ 
&  &  &  &  &  &  & \multicolumn{1}{r}{} \\ 
LQML & bias & \multicolumn{1}{r}{$-0.0004$} & \multicolumn{1}{r}{$-0.0008$}
& \multicolumn{1}{r}{$-0.0004$} & \multicolumn{1}{r}{$-0.0006$} & 
\multicolumn{1}{r}{$-0.0003$} & \multicolumn{1}{r}{$-0.0001$} \\ 
& \textit{rmse} & \multicolumn{1}{r}{$\mathit{0.0248}$} & \multicolumn{1}{r}{%
$\mathit{0.0244}$} & \multicolumn{1}{r}{$\mathit{0.0229}$} & 
\multicolumn{1}{r}{$\mathit{0.0184}$} & \multicolumn{1}{r}{$\mathit{0.0102}$}
& \multicolumn{1}{r}{$\mathit{0.0066}$} \\ 
&  &  &  &  &  &  & \multicolumn{1}{r}{} \\ 
DQML$_{\boldsymbol{x}}$ & bias & \multicolumn{1}{r}{$-0.0030$} & 
\multicolumn{1}{r}{$-0.0048$} & \multicolumn{1}{r}{$-0.0060$} & 
\multicolumn{1}{r}{$-0.0088$} & \multicolumn{1}{r}{$-0.0083$} & 
\multicolumn{1}{r}{$-0.0063$} \\ 
& \textit{rmse} & \multicolumn{1}{r}{$\mathit{0.0328}$} & \multicolumn{1}{r}{%
$\mathit{0.0320}$} & \multicolumn{1}{r}{$\mathit{0.0299}$} & 
\multicolumn{1}{r}{$\mathit{0.0244}$} & \multicolumn{1}{r}{$\mathit{0.0152}$}
& \multicolumn{1}{r}{$\mathit{0.0103}$} \\ 
& \multicolumn{1}{c}{} & \multicolumn{1}{r}{} & \multicolumn{1}{r}{} & 
\multicolumn{1}{r}{} &  & \multicolumn{1}{r}{} & \multicolumn{1}{r}{} \\ 
DQML$_{\Delta \boldsymbol{x}}$ & bias & \multicolumn{1}{r}{$-0.0030$} & 
\multicolumn{1}{r}{$-0.0048$} & \multicolumn{1}{r}{$-0.0060$} & 
\multicolumn{1}{r}{$-0.0088$} & \multicolumn{1}{r}{$-0.0083$} & 
\multicolumn{1}{r}{$-0.0063$} \\ 
& \textit{rmse} & \multicolumn{1}{r}{$\mathit{0.0328}$} & \multicolumn{1}{r}{%
$\mathit{0.0320}$} & \multicolumn{1}{r}{$\mathit{0.0299}$} & 
\multicolumn{1}{r}{$\mathit{0.0244}$} & \multicolumn{1}{r}{$\mathit{0.0152}$}
& \multicolumn{1}{r}{$\mathit{0.0103}$} \\ 
&  &  &  &  &  &  & \multicolumn{1}{r}{} \\ \hline

\end{tabular}
\end{center}

\vspace{.75cm}

In order for the system GMM estimator to be consistent (as $N\rightarrow
\infty $) the stochastic process for each individual has to have had
sufficient time to converge to its steady state by time $t=1$ (see, e.g.,
Roodman 2009). However, given $t_{0}=1$, convergence to a steady state at
time $t=1$ has clearly not occurred. The effect of the failure of this
initial condition restriction is most striking when $\sigma _{\zeta }=4$. In this case, for many designs, the absolute bias and root mean squared error of the system GMM estimator is much larger than that of the other estimators.

Except for the condition that $y_{i0}$ must be uncorrelated with $v_{it}$
for $t\geq 1$, the QML estimators are unaffected by initial conditions.
However, the consistency (as $N\rightarrow \infty $) of the differenced QML
estimators --- DQML$_{\boldsymbol{x}}$ and DQML$_{\Delta \boldsymbol{x}}$
--- depends on the $v_{it}$s being unconditionally homoskedastic, and, when $%
t_{0}=1$, the $v_{it}$s are time-series heteroskedastic. Consequently, in
Table 2, the differenced QML estimators no longer dominate the
differenced GMM estimator in terms of finite sample bias. On the other hand,
the levels QML\ estimator is robust with respect to time-series
heteroskedasticity, and therefore its finite sample bias is still negligible
for $t_{0}=1$.

\section{Conclusions}

This paper established the almost sure convergence and asymptotic normality
of levels and differenced QML estimators of the parameters of a $p$th-order
dynamic panel data model. The almost sure convergence and asymptotic
normality of the estimators do not depend on initial conditions, like those
required by the sytem GMM estimator. Moreover, the log-likelihood can be
misspecified, and the errors can be conditionally and time-series
heteroskedastic. However, only levels QML estimates can be easily calculated
when the errors are time-series heteroskedastic. The paper provided an
ECME algorithm for this case. Furthermore, the levels QML estimator
dominated all of the other estimators in terms of having the smallest root
mean squared errors.

\newpage
\section*{Appendix A: Lemma 1 Proof}

In order to establish $E\left( \boldsymbol{y}_{i,-j}^{\prime
}\Omega _{0}^{\ast -1}\boldsymbol{e}_{i}\right) =0$, I first use an analysis
similar to that in Hamilton (1994, pp. 7-9). Let $\boldsymbol{\xi }%
_{it}=\left( y_{it},y_{i,t-1},\ldots ,y_{i,t-p+1}\right) ^{\prime }$, $%
\boldsymbol{\varsigma }_{it}=\left( \boldsymbol{x}_{it}^{\prime }\boldsymbol{%
\beta }_{0}+e_{it},0,\ldots ,0\right) ^{\prime }$, and%
\begin{equation}
\boldsymbol{F}=\left( 
\begin{array}{ccccc}
\delta _{01} & \delta _{02} & \cdots & \delta _{0,p-1} & \delta _{0p} \\ 
1 & 0 & \cdots & 0 & 0 \\ 
0 & 1 & \cdots & 0 & 0 \\ 
\vdots & \vdots & \ddots & \vdots & \vdots \\ 
0 & 0 & \cdots & 1 & 0%
\end{array}%
\right) ,  \label{F_def}
\end{equation}%
where $\boldsymbol{\delta }_{0}=\left( \delta _{01},\ldots ,\delta
_{0p}\right) ^{\prime }$. Then $\boldsymbol{\xi }_{it}=\boldsymbol{F\xi }%
_{i,t-1}+\boldsymbol{\varsigma }_{it}$. Hence, $\boldsymbol{\xi }_{i1}=%
\boldsymbol{F\xi }_{i0}+\boldsymbol{\varsigma }_{i1}$, and, for $t>1$, by
repeated substitutions we get $\boldsymbol{\xi }_{it}=\boldsymbol{F}^{t}%
\boldsymbol{\xi }_{i0}+\boldsymbol{F}^{t-1}\boldsymbol{\varsigma }_{i1}+%
\boldsymbol{F}^{t-2}\boldsymbol{\varsigma }_{i2}+\cdots +\boldsymbol{%
F\varsigma }_{i,t-1}+\boldsymbol{\varsigma }_{it}$. Writing this last
expression out in full, we have%
\begin{eqnarray}
\left( 
\begin{array}{c}
y_{it} \\ 
y_{i,t-1} \\ 
\vdots \\ 
y_{i,t-p+1}%
\end{array}%
\right) &=&\boldsymbol{F}^{t}\left( 
\begin{array}{c}
y_{i0} \\ 
y_{i,-1} \\ 
\vdots \\ 
y_{i,-p+1}%
\end{array}%
\right) +\boldsymbol{F}^{t-1}\left( 
\begin{array}{c}
\boldsymbol{x}_{i1}^{\prime }\boldsymbol{\beta }_{0}+e_{i1} \\ 
0 \\ 
\vdots \\ 
0%
\end{array}%
\right) +\boldsymbol{F}^{t-2}\left( 
\begin{array}{c}
\boldsymbol{x}_{i2}^{\prime }\boldsymbol{\beta }_{0}+e_{i2} \\ 
0 \\ 
\vdots \\ 
0%
\end{array}%
\right)  \notag \\
&&+\cdots +\boldsymbol{F}\left( 
\begin{array}{c}
\boldsymbol{x}_{i,t-1}^{\prime }\boldsymbol{\beta }_{0}+e_{i,t-1} \\ 
0 \\ 
\vdots \\ 
0%
\end{array}%
\right) +\left( 
\begin{array}{c}
\boldsymbol{x}_{it}^{\prime }\boldsymbol{\beta }_{0}+e_{it} \\ 
0 \\ 
\vdots \\ 
0%
\end{array}%
\right) .  \label{xi_it_full}
\end{eqnarray}%
Next let $f_{rs}^{\left( t\right) }$ denote the $\left( r,s\right) $th
element of $\boldsymbol{F}^{t}$. Then $y_{i1}=f_{11}^{\left( 1\right)
}y_{i0}+f_{12}^{\left( 1\right) }y_{i,-1}+\cdots +f_{1p}^{\left( 1\right)
}y_{i,-p+1}+\boldsymbol{x}_{i1}^{\prime }\boldsymbol{\beta }_{0}+e_{i1}$,
and, for $t>1$, from the first equation in (\ref{xi_it_full}) we see that%
\begin{eqnarray}
y_{it} &=&f_{11}^{\left( t\right) }y_{i0}+f_{12}^{\left( t\right)
}y_{i,-1}+\cdots +f_{1p}^{\left( t\right) }y_{i,-p+1}+f_{11}^{\left(
t-1\right) }\left( \boldsymbol{x}_{i1}^{\prime }\boldsymbol{\beta }%
_{0}+e_{i1}\right)  \notag \\
&&+f_{11}^{\left( t-2\right) }\left( \boldsymbol{x}_{i2}^{\prime }%
\boldsymbol{\beta }_{0}+e_{i2}\right) +\cdots +f_{11}^{\left( 1\right)
}\left( \boldsymbol{x}_{i,t-1}^{\prime }\boldsymbol{\beta }%
_{0}+e_{i,t-1}\right) +\boldsymbol{x}_{it}^{\prime }\boldsymbol{\beta }%
_{0}+e_{it}.  \label{y_it_0}
\end{eqnarray}

Using the expression for $y_{it}$ in Eq. (\ref{y_it_0}), we can write $%
\boldsymbol{y}_{i,-j}$ in terms of $\boldsymbol{y}_{i}^{o}$, $\boldsymbol{X}%
_{i}$, and $\boldsymbol{e}_{i}$. To that end, let $\boldsymbol{A}_{j}$ and $%
\boldsymbol{B}_{j}$ be $T\times p$ and $T\times T$ matrices given by%
\begin{equation}
\boldsymbol{A}_{j}=\left( 
\begin{array}{cccccccc}
0 & 0 & \cdots & 0 & 1 & 0 & \cdots & 0 \\ 
0 & 0 & \cdots & 1 & 0 & 0 & \cdots & 0 \\ 
\vdots & \vdots &  & \vdots & \vdots & \vdots &  & \vdots \\ 
0 & 1 & \cdots & 0 & 0 & 0 & \cdots & 0 \\ 
1 & 0 & \cdots & 0 & 0 & 0 & \cdots & 0 \\ 
f_{11}^{\left( 1\right) } & f_{12}^{\left( 1\right) } & \cdots & 
f_{1,j-1}^{\left( 1\right) } & f_{1j}^{\left( 1\right) } & f_{1,j+1}^{\left(
1\right) } & \cdots & f_{1p}^{\left( 1\right) } \\ 
f_{11}^{\left( 2\right) } & f_{12}^{\left( 2\right) } & \cdots & 
f_{1,j-1}^{\left( 2\right) } & f_{1j}^{\left( 2\right) } & f_{1,j+1}^{\left(
2\right) } & \cdots & f_{1p}^{\left( 2\right) } \\ 
\vdots & \vdots &  & \vdots & \vdots & \vdots &  & \vdots \\ 
f_{11}^{\left( T-j\right) } & f_{12}^{\left( T-j\right) } & \cdots & 
f_{1,j-1}^{\left( T-j\right) } & f_{1j}^{\left( T-j\right) } & 
f_{1,j+1}^{\left( T-j\right) } & \cdots & f_{1p}^{\left( T-j\right) }%
\end{array}%
\right)  \label{A_j}
\end{equation}%
\begin{equation}
\boldsymbol{B}_{j}=\left( 
\begin{array}{ccccccccc}
0 & 0 & 0 & \cdots & 0 & 0 & 0 & \cdots & 0 \\ 
\vdots & \vdots & \vdots &  & \vdots & \vdots & \vdots &  & \vdots \\ 
0 & 0 & 0 & \cdots & 0 & 0 & 0 & \cdots & 0 \\ 
1 & 0 & 0 & \cdots & 0 & 0 & 0 & \cdots & 0 \\ 
f_{11}^{\left( 1\right) } & 1 & 0 & \cdots & 0 & 0 & 0 & \cdots & 0 \\ 
f_{11}^{\left( 2\right) } & f_{11}^{\left( 1\right) } & 1 & \cdots & 0 & 0 & 
0 & \cdots & 0 \\ 
\vdots & \vdots & \vdots &  & \vdots & \vdots & \vdots &  & \vdots \\ 
f_{11}^{\left( T-j-1\right) } & f_{11}^{\left( T-j-2\right) } & 
f_{11}^{\left( T-j-3\right) } & \cdots & f_{11}^{\left( 1\right) } & 1 & 0 & 
\cdots & 0%
\end{array}%
\right) .  \label{B_j}
\end{equation}%
Given these definitions, we have $\boldsymbol{y}_{i,-j}=\boldsymbol{A}_{j}%
\boldsymbol{y}_{i}^{o}+\boldsymbol{B}_{j}\left( \boldsymbol{X}_{i}%
\boldsymbol{\beta }_{0}+\boldsymbol{e}_{i}\right) $.

Therefore, $E\left( \boldsymbol{y}_{i,-j}^{\prime }\Omega _{0}^{\ast -1}%
\boldsymbol{e}_{i}\right) =E\left( \boldsymbol{y}_{i}^{o\prime }\boldsymbol{A%
}_{j}^{\prime }\Omega _{0}^{\ast -1}\boldsymbol{e}_{i}\right) +E\left( 
\boldsymbol{\beta }_{0}^{\prime }\boldsymbol{X}_{i}^{\prime }\boldsymbol{B}%
_{j}^{\prime }\Omega _{0}^{\ast -1}\boldsymbol{e}_{i}\right) +E\left( 
\boldsymbol{e}_{i}^{\prime }\boldsymbol{B}_{j}^{\prime }\Omega _{0}^{\ast -1}%
\boldsymbol{e}_{i}\right) $. Note that $E\left( \boldsymbol{e}_{i}^{\prime }%
\boldsymbol{B}_{j}^{\prime }\Omega _{0}^{\ast -1}\boldsymbol{e}_{i}\right) =E%
\left[ \text{tr}\left( \Omega _{0}^{\ast -1}\boldsymbol{e}_{i}\boldsymbol{e}%
_{i}^{\prime }\boldsymbol{B}_{j}^{\prime }\right) \right] =$ tr$\left[
\Omega _{0}^{\ast -1}E\left( \boldsymbol{e}_{i}\boldsymbol{e}_{i}^{\prime
}\right) \boldsymbol{B}_{j}^{\prime }\right] =$ tr$\left( \boldsymbol{B}%
_{j}^{\prime }\right) =0$, where the last equality follows from the fact
that $\boldsymbol{B}_{j}$ is a square matrix with zeros down the main
diagonal. Moreover, if $E\left( \boldsymbol{e}_{i}\boldsymbol{x}_{i}^{\prime
}\right) =\mathbf{0}$, then $E\left( \boldsymbol{\beta }_{0}^{\prime }%
\boldsymbol{X}_{i}^{\prime }\boldsymbol{B}_{j}^{\prime }\Omega _{0}^{\ast -1}%
\boldsymbol{e}_{i}\right) =0$. And $E\left( \boldsymbol{y}_{i}^{o\prime }%
\boldsymbol{A}_{j}^{\prime }\Omega _{0}^{\ast -1}\boldsymbol{e}_{i}\right) =$
tr$\left[ E\left( \boldsymbol{e}_{i}\boldsymbol{y}_{i}^{o\prime }\right) 
\boldsymbol{A}_{j}^{\prime }\Omega _{0}^{\ast -1}\right] =0$ given $E\left( 
\boldsymbol{e}_{i}\boldsymbol{y}_{i}^{o\prime }\right) =\mathbf{0}$. The
preceding proves $E\left( \boldsymbol{y}_{i,-j}^{\prime }\Omega _{0}^{\ast
-1}\boldsymbol{e}_{i}\right) =0$.%

\section*{Appendix B: Theorem 1 Proof}

The proof of Theorem 1 relies on verifying several preliminary
results, which are provided as Lemmas B.1 through B.3. Throughout
convergence is with respect to $N\rightarrow \infty $, with $T$ fixed.
Moreover, in the sequel, $M$ denotes a sufficiently large finite number.

\bigskip

\noindent \textbf{Lemma B.1}. Suppose $E\left( x_{itk}^{2}\right) <\infty $
and $E\left( y_{it}^{2}\right) <\infty $, for each $i$, $t$, and $k$, and
Conditions C2 and C4 are satisfied. Then the linear projection in (\ref{lp})
exists. Furthermore, the limits $L\left( \boldsymbol{\psi }\right)
=\lim_{N\rightarrow \infty }E\left[ L_{N}\left( \boldsymbol{\psi }\right) %
\right] $ and $\boldsymbol{H}\left( \boldsymbol{\psi }\right)
=\lim_{N\rightarrow \infty }E\left[ \boldsymbol{H}_{N}\left( \boldsymbol{%
\psi }\right) \right] $ exist, and $L\left( \boldsymbol{\psi }\right) $ and
the elements of $\boldsymbol{H}\left( \boldsymbol{\psi }\right) $ are
continuous functions of $\boldsymbol{\psi }$.

\bigskip

\noindent \textbf{Proof}. The conditions $E\left( x_{itk}^{2}\right) <\infty 
$ and $E\left( y_{it}^{2}\right) <\infty $, for each $i$, $t$, and $k$, and
C2 imply the existence of the linear projection in (\ref{lp}) (see, e.g.,
Wooldridge, 2010, pp. 25-26).%

Also, $E\left[ L_{N}\left( \boldsymbol{\psi }\right) \right] $ is
finite if $E\left[ \boldsymbol{u}_{i}\left( \boldsymbol{\gamma }\right)
^{\prime }\Omega ^{-1}\boldsymbol{u}_{i}\left( \boldsymbol{\gamma }\right) %
\right] $ is finite, and the latter is finite if $x_{itk}$ and $y_{it}$ have
finite second-order moments, for all $i$, $t$, and $k$.

\sloppy The matrix $E\left[ \boldsymbol{H}_{N}\left( \boldsymbol{\psi }\right) %
\right] $ has finite elements as well. To see this, first let $\boldsymbol{W}%
_{i\cdot j}$ denote the $j$th column of $\boldsymbol{W}_{i}$, and let $%
\boldsymbol{S}_{\cdot j}$ denote the $j$th column of $\partial $vec$\left(
\Omega \right) /\partial \boldsymbol{\omega }^{\prime }$, where recall that $%
\boldsymbol{\omega }=$ vech$\left( \Omega \right) $. Then, $\partial
^{2}l_{i}\left( \boldsymbol{\psi }\right) /\partial \gamma _{j}\partial
\gamma _{k}=-\boldsymbol{W}_{i\cdot j}^{\prime }\Omega ^{-1}\boldsymbol{W}%
_{i\cdot k}$, $\partial ^{2}l_{i}\left( \boldsymbol{\psi }\right) /\partial
\gamma _{j}\partial \omega _{k}=-\boldsymbol{W}_{i\cdot j}^{\prime }\Omega
^{-1}\left( \partial \Omega /\partial \omega _{k}\right) \Omega ^{-1}%
\boldsymbol{u}_{i}\left( \boldsymbol{\gamma }\right) $, and%
\begin{equation}
\frac{\partial ^{2}l_{i}\left( \boldsymbol{\psi }\right) }{\partial \omega
_{j}\partial \omega _{k}}=\frac{1}{2}\boldsymbol{S}_{\cdot j}^{\prime
}\left( \Omega ^{-1}\otimes \Omega ^{-1}\right) \boldsymbol{S}_{\cdot k}-%
\frac{1}{2}s_{ijk}^{\left( 1\right) }\left( \boldsymbol{\psi }\right) -\frac{%
1}{2}s_{ijk}^{\left( 2\right) }\left( \boldsymbol{\psi }\right) ,
\label{2nd_deriv_3}
\end{equation}%
where $s_{ijk}^{\left( 1\right) }\left( \boldsymbol{\psi }\right) =%
\boldsymbol{S}_{\cdot j}^{\prime }\left( \Omega ^{-1}\otimes \Omega ^{-1}%
\boldsymbol{u}_{i}\left( \boldsymbol{\gamma }\right) \boldsymbol{u}%
_{i}\left( \boldsymbol{\gamma }\right) ^{\prime }\Omega ^{-1}\right) 
\boldsymbol{S}_{\cdot k}$ and $s_{ijk}^{\left( 2\right) }\left( \boldsymbol{%
\psi }\right) =\boldsymbol{S}_{\cdot j}^{\prime }\left( \Omega ^{-1}%
\boldsymbol{u}_{i}\left( \boldsymbol{\gamma }\right) \boldsymbol{u}%
_{i}\left( \boldsymbol{\gamma }\right) ^{\prime }\Omega ^{-1}\otimes \Omega
^{-1}\right) \boldsymbol{S}_{\cdot k}$ (see Ruud 2000, p. 930). From the
preceding second-order partial derivatives we see that the condition $%
E\left( x_{itk}^{2}\right) <\infty $ and $E\left( y_{it}^{2}\right) <\infty $%
, for each $i$, $t$, and $k$, implies $E\left[ \boldsymbol{H}_{N}\left( 
\boldsymbol{\psi }\right) \right] $ has finite elements.%

Inspection of $E\left[ L_{N}\left( \boldsymbol{\psi }\right) \right] $ and
the elements of $E\left[ \boldsymbol{H}_{N}\left( \boldsymbol{\psi }\right) %
\right] $ reveals $E\left[ L_{N}\left( \boldsymbol{\psi }\right) \right] $
and the elements of $E\left[ \boldsymbol{H}_{N}\left( \boldsymbol{\psi }%
\right) \right] $ are functions of $\boldsymbol{\psi }$ and terms of the
form $N^{-1}\sum_{i}E\left( y_{is}y_{it}\right) $, $N^{-1}\sum_{i}E\left(
y_{is}x_{itk}\right) $, and $N^{-1}\sum_{i}E\left( x_{isj}x_{itk}\right) $.
Therefore, if the limits of these averages exist (as $N\rightarrow \infty $%
), then the limits $L\left( \boldsymbol{\psi }\right) =\lim_{N\rightarrow
\infty }E\left[ L_{N}\left( \boldsymbol{\psi }\right) \right] $ and $%
\boldsymbol{H}\left( \boldsymbol{\psi }\right) =\lim_{N\rightarrow \infty }E%
\left[ \boldsymbol{H}_{N}\left( \boldsymbol{\psi }\right) \right] $ exist,
where $L\left( \boldsymbol{\psi }\right) $ and the elements of $\boldsymbol{H%
}\left( \boldsymbol{\psi }\right) $ are functions of $\boldsymbol{\psi }$
and terms involving limits of the form $\lim_{N\rightarrow \infty
}N^{-1}\sum_{i}E\left( y_{is}y_{it}\right) $, $\lim_{N\rightarrow \infty
}N^{-1}\sum_{i}E\left( y_{is}x_{itk}\right) $, and $\lim_{N\rightarrow
\infty }N^{-1}\sum_{i}E\left( x_{isj}x_{itk}\right) $. And, inspection of $%
L\left( \boldsymbol{\psi }\right) $ and the elements of $\boldsymbol{H}%
\left( \boldsymbol{\psi }\right) $ reveals $L\left( \boldsymbol{\psi }%
\right) $ and the elements of $\boldsymbol{H}\left( \boldsymbol{\psi }%
\right) $ are continuous functions of $\boldsymbol{\psi }$.%

\bigskip
\fussy
\noindent \textbf{Lemma B.2}. Let $\overline{\Psi }$ denote a compact subset
of $\Psi $. Suppose C1, C4, and C5 are satisfied. Then $L_{N}\left( 
\boldsymbol{\cdot }\right) \overset{a.s.}{\rightarrow }L\left( \boldsymbol{%
\cdot }\right) $ uniformly on $\overline{\Psi }$.

\bigskip
\sloppy
\noindent \textbf{Proof}. Let $\omega ^{st}$ denote the $\left( s,t\right) $%
th element of $\Omega ^{-1}$; let $\gamma _{k}$ denote the $k$th element of $%
\boldsymbol{\gamma }$; recall that $\boldsymbol{W}_{i\cdot j}$ is the $j$th
column of $\boldsymbol{W}_{i}$; and let $W_{itj}$ denote the $t$th element
of $\boldsymbol{W}_{i\cdot j}$. Also, let $S_{y_{s}y_{t},N}=N^{-1}\sum_{i}%
\left[ y_{is}y_{it}-E\left( y_{is}y_{it}\right) \right] $, $%
S_{y_{s}W_{tj},N}=N^{-1}\sum_{i}\left[ y_{is}W_{itj}-E\left(
y_{is}W_{itj}\right) \right] $, and $S_{W_{sj}W_{tk},N}=N^{-1}\sum_{i}\left[
W_{isj}W_{itk}-E\left( W_{isj}W_{itk}\right) \right] $. Then $L_{N}\left( 
\boldsymbol{\psi }\right) -E\left[ L_{N}\left( \boldsymbol{\psi }\right) %
\right] =-\sum_{s}\sum_{t}\omega
^{st}S_{y_{s}y_{t},N}/2+\sum_{s}\sum_{t}\omega ^{st}\sum_{j}\gamma
_{j}S_{y_{s}W_{tj},N}-\sum_{s}\sum_{t}\omega ^{st}\sum_{j}\sum_{k}\gamma
_{j}\gamma _{k}S_{W_{sj}W_{tk},N}/2$. Therefore, by an obvious inequality,
we have $\left\vert L_{N}\left( \boldsymbol{\psi }\right) -E\left[
L_{N}\left( \boldsymbol{\psi }\right) \right] \right\vert \leq
\sum_{s}\sum_{t}\left\vert \omega ^{st}\right\vert \left\vert
S_{y_{s}y_{t},N}\right\vert /2+\sum_{s}\sum_{t}\left\vert \omega
^{st}\right\vert \sum_{k}\left\vert \gamma _{k}\right\vert \left\vert
S_{y_{s}W_{tk},N}\right\vert +\sum_{s}\sum_{t}\left\vert \omega
^{st}\right\vert \sum_{j}\sum_{k}\left\vert \gamma _{j}\gamma
_{k}\right\vert \left\vert S_{W_{sj}W_{tk},N}\right\vert /2$. Given $\omega
^{st}$ and $\gamma _{k}$ are bounded for $\boldsymbol{\psi }\in \overline{%
\Psi }$, it follows that%
\begin{eqnarray}
\sup_{\boldsymbol{\psi }\in \overline{\Psi }}\left\vert L_{N}\left( 
\boldsymbol{\psi }\right) -E\left[ L_{N}\left( \boldsymbol{\psi }\right) %
\right] \right\vert &\leq &M\sum_{s}\sum_{t}\left\vert
S_{y_{s}y_{t},N}\right\vert +M\sum_{s}\sum_{t}\sum_{k}\left\vert
S_{y_{s}W_{tk},N}\right\vert  \notag \\
&&+M\sum_{s}\sum_{t}\sum_{j}\sum_{k}\left\vert S_{W_{sj}W_{tk},N}\right\vert
.  \label{unif_conv_pr_1}
\end{eqnarray}%
Hence, $L_{N}\left( \cdot \right) -E\left[ L_{N}\left( \cdot \right) \right] 
\overset{a.s.}{\rightarrow }0$ uniformly on $\overline{\Psi }$ if $%
S_{y_{s}y_{t},N}\overset{a.s.}{\rightarrow }0$, $S_{y_{s}W_{tk},N}\overset{%
a.s.}{\rightarrow }0$, and $S_{W_{sj}W_{tk},N}\overset{a.s.}{\rightarrow }0$
for each $s$, $t$, $j$, and $k$.

\fussy To see that $S_{y_{s}y_{t},N}\overset{a.s.}{\rightarrow }0$, note that, by
the Cauchy-Schwarz inequality and C1, we get $E\left\vert
y_{is}y_{it}\right\vert ^{1+\epsilon /2}\leq \left( E\left\vert
y_{is}\right\vert ^{2+\epsilon }E\left\vert y_{it}\right\vert ^{2+\epsilon
}\right) ^{1/2}<M$ for some $\epsilon >0$ and all $i$, $s$, and $t$. This
conclusion and C5 imply $S_{y_{s}y_{t},N}\overset{a.s.}{\rightarrow }0$ (see
White 2001, p. 35, Corollary 3.9). By similar arguments, we also have $%
S_{y_{s}W_{tk},N}\overset{a.s.}{\rightarrow }0$ and $S_{W_{sj}W_{tk},N}%
\overset{a.s.}{\rightarrow }0$. Hence, $L_{N}\left( \cdot \right) -E\left[
L_{N}\left( \cdot \right) \right] \overset{a.s.}{\rightarrow }0$ uniformly
on $\overline{\Psi }$.

\sloppy Given C4, the following expressions are defined: $A_{y_{s}y_{t},N}=N^{-1}%
\sum_{i}E\left( y_{is}y_{it}\right) -\lim_{N\rightarrow \infty
}N^{-1}\sum_{i}E\left( y_{is}y_{it}\right) $, $A_{y_{s}W_{tj},N}=N^{-1}%
\sum_{i}E\left( y_{is}W_{itj}\right) -\lim_{N\rightarrow \infty
}N^{-1}\sum_{i}E\left( y_{is}W_{itj}\right) $, and $%
A_{W_{sj}W_{tk},N}=N^{-1}\sum_{i}E\left( W_{isj}W_{itk}\right)
-\lim_{N\rightarrow \infty }N^{-1}\sum_{i}E\left( W_{isj}W_{itk}\right) $.
And, by arguments analogous to those leading to the inequality in (\ref%
{unif_conv_pr_1}), one can show $\sup_{\boldsymbol{\psi }\in \overline{\Psi }%
}\left\vert E\left[ L_{N}\left( \boldsymbol{\psi }\right) \right] -L\left( 
\boldsymbol{\psi }\right) \right\vert \leq M\sum_{s}\sum_{t}\left\vert
A_{y_{s}y_{t},N}\right\vert +M\sum_{s}\sum_{t}\sum_{j}\left\vert
A_{y_{s}W_{tj},N}\right\vert +M\sum_{s}\sum_{t}\sum_{j}\sum_{k}\left\vert
A_{W_{sj}W_{tk},N}\right\vert $. Because $A_{y_{s}y_{t},N}$, $%
A_{y_{s}W_{tj},N}$, and $A_{W_{sj}W_{tk},N}$ all $\rightarrow 0$, we have $E%
\left[ L_{N}\left( \cdot \right) \right] \rightarrow L\left( \cdot \right) $
uniformly on $\overline{\Psi }$.

\fussy The conclusions of the last two paragraphs imply $L_{N}\left( \cdot
\right) \overset{a.s.}{\rightarrow }L\left( \cdot \right) $ uniformly on $%
\overline{\Psi }$.%

\bigskip

\noindent \textbf{Lemma B.3}. If C1--C3 are satisfied, then $E\left[\partial L_{N}\left( \boldsymbol{\psi }_{0}\right) /\partial 
\boldsymbol{\psi }\right] = \bold{0}$. If, in addition, C4 and C5 are satisfied and $%
\boldsymbol{H}_{0}$ is negative definite, then there is a compact subset $%
\overline{\Psi }$ of $\Psi $, with $\boldsymbol{\psi }_{0}$ in its interior,
such that $L\left( \boldsymbol{\psi }\right) <L\left( \boldsymbol{\psi }%
_{0}\right) $ if $\boldsymbol{\psi }\in \overline{\Psi }$ and $\boldsymbol{%
\psi }\neq \boldsymbol{\psi }_{0}$.

\bigskip

\noindent \textbf{Proof}. First%
\begin{equation}
E\left[ \partial l_{i}\left( \boldsymbol{\psi }_{0}\right) /\partial 
\boldsymbol{\psi }\right] =\mathbf{0}  \label{gradient1}
\end{equation}%
is established. By well known results, $\partial l_{i}\left( \boldsymbol{%
\psi }\right) /\partial \boldsymbol{\gamma }=\boldsymbol{W}_{i}^{\prime
}\Omega ^{-1}\boldsymbol{u}_{i}\left( \boldsymbol{\gamma }\right) $ and%
\begin{equation}
\frac{\partial l_{i}\left( \boldsymbol{\psi }\right) }{\partial \boldsymbol{%
\omega }}=-\frac{1}{2}\text{vech}\left( \Omega ^{-1}-\Omega ^{-1}\boldsymbol{%
u}_{i}\left( \boldsymbol{\gamma }\right) \boldsymbol{u}_{i}\left( 
\boldsymbol{\gamma }\right) ^{\prime }\Omega ^{-1}\right)  \label{gradient2}
\end{equation}%
(see, e.g., Ruud, 2000, pp. 928-930). To see that $E\left[ \partial
l_{i}\left( \boldsymbol{\psi }_{0}\right) /\partial \boldsymbol{\gamma }%
\right] =\mathbf{0}$, first note that $E\left( \boldsymbol{Z}_{i}^{\prime
}\Omega _{0}^{-1}\boldsymbol{u}_{i}\right) =\mathbf{0}$ because all of the
elements of $\boldsymbol{u}_{i}$ are uncorrelated with all of the elements
of $\boldsymbol{Z}_{i}$ by construction. Moreover, C1--C3 imply $%
E\left( \boldsymbol{u}_{i}\boldsymbol{y}_{i}^{o\prime }\right) =\mathbf{0}$
and $E\left( \boldsymbol{u}_{i}\boldsymbol{x}_{i}^{\prime }\right) =\mathbf{0%
}$, and $E\left( \boldsymbol{u}_{i}\boldsymbol{u}_{i}^{\prime }\right)
=\Omega _{0}$. Thus, the conditions of Lemma 1 hold for the augmented
regression in (\ref{aug_model}). Hence, by Lemma 1, we have $E\left( 
\boldsymbol{y}_{i,-j}^{\prime }\Omega _{0}^{-1}\boldsymbol{u}_{i}\right) =0$
($j=1,\ldots ,p$). This proves $E\left[ \partial l_{i}\left( \boldsymbol{%
\psi }\right) /\partial \boldsymbol{\gamma }\right] =\mathbf{0}$.
Furthermore, from Eq. (\ref{gradient2}), it is clear that, because $E\left( 
\boldsymbol{u}_{i}\boldsymbol{u}_{i}^{\prime }\right) =\Omega _{0}$, we have 
$E\left[ \partial l_{i}\left( \boldsymbol{\psi }_{0}\right) /\partial 
\boldsymbol{\omega }\right] =\mathbf{0}$. Hence, $E\left[\partial L_{N}\left( \boldsymbol{\psi }_{0}\right) /\partial 
\boldsymbol{\psi }\right] = \bold{0}$.

Next, a Taylor series expansion gives%
\begin{equation}
L_{N}\left( \boldsymbol{\psi }\right) =L_{N}\left( \boldsymbol{\psi }%
_{0}\right) +\left( \boldsymbol{\psi }-\boldsymbol{\psi }_{0}\right)
^{\prime }\boldsymbol{g}_{N}\left( \boldsymbol{\psi }_{0}\right) +\left( 
\boldsymbol{\psi }-\boldsymbol{\psi }_{0}\right) ^{\prime }\boldsymbol{H}%
_{N}\left( \boldsymbol{\psi }^{\ast }\right) \left( \boldsymbol{\psi }-%
\boldsymbol{\psi }_{0}\right) /2,  \label{Taylor0}
\end{equation}%
where $\boldsymbol{g}_{N}\left( \boldsymbol{\psi }_{0}\right) =\partial
L_{N}\left( \boldsymbol{\psi }_{0}\right) /\partial \boldsymbol{\psi }$, and 
$\boldsymbol{\psi }^{\ast }$ satisfies $\left\Vert \boldsymbol{\psi }-%
\boldsymbol{\psi }^{\ast }\right\Vert \leq \left\Vert \boldsymbol{\psi }-%
\boldsymbol{\psi }_{0}\right\Vert $. Given Eq. (\ref{gradient1}) and Lemma
B.1, taking the expectation of the left and right-hand sides of (\ref%
{Taylor0}) and then letting $N\rightarrow \infty $ gives $L\left( 
\boldsymbol{\psi }\right) =L\left( \boldsymbol{\psi }_{0}\right) +\left( 
\boldsymbol{\psi }-\boldsymbol{\psi }_{0}\right) ^{\prime }\boldsymbol{H}%
\left( \boldsymbol{\psi }^{\ast }\right) \left( \boldsymbol{\psi }-%
\boldsymbol{\psi }_{0}\right) /2$.

Let $h_{jk}\left( \boldsymbol{\psi }\right) $ denote the $\left( j,k\right) $%
th element of $\boldsymbol{H}\left( \boldsymbol{\psi }\right) $, and define
determinants%
\begin{equation*}
d_{j}\left( \boldsymbol{\psi }\right) =\left\vert 
\begin{array}{ccc}
h_{11}\left( \boldsymbol{\psi }\right) & \cdots & h_{1j}\left( \boldsymbol{%
\psi }\right) \\ 
\vdots & \ddots & \vdots \\ 
h_{j1}\left( \boldsymbol{\psi }\right) & \cdots & h_{jj}\left( \boldsymbol{%
\psi }\right)%
\end{array}%
\right\vert \text{ \ \ \ \ \ \ \ \ \ \ \ \ \ }\left( j=1,\ldots ,m\right) .
\end{equation*}%
By assumption, $%
\boldsymbol{H}_{0}=\boldsymbol{H}\left( \boldsymbol{\psi }_{0}\right) $ is
negative definite, and thus $d_{1}\left( \boldsymbol{\psi }_{0}\right) <0$, $%
d_{2}\left( \boldsymbol{\psi }_{0}\right) >0$, $d_{3}\left( \boldsymbol{\psi 
}_{0}\right) <0,\ldots $ (see Rao 1973, p. 37). Moreover, the determinant $d_{j}\left( \cdot
\right) $ is continuous in $h_{11}\left( \cdot \right) ,$ $h_{12}\left(
\cdot \right) ,\ldots $, which are, in turn, continuous in $\boldsymbol{\psi 
}$ (see Lemma B.1). Hence, $d_{j}\left( \cdot \right) $ is continuous in $%
\boldsymbol{\psi }$. It follows that there is a $r>0$ such that for the
closed ball in $\boldsymbol{%
\mathbb{R}
}^{m}$, centered at $\boldsymbol{\psi }_{0}$, with radius $r$, we have $%
d_{1}\left( \boldsymbol{\psi }\right) <0$, $d_{2}\left( \boldsymbol{\psi }%
\right) >0$, $d_{3}\left( \boldsymbol{\psi }\right) <0,\ldots $ for $%
\boldsymbol{\psi }$ in the ball. Let $\overline{\Psi }$ denote the ball (a
compact subset of $%
\mathbb{R}
^{m}$). Then $\boldsymbol{H}\left( \boldsymbol{\psi }\right) $ is negative
definite for $\boldsymbol{\psi }\in \overline{\Psi }$. Therefore, for $%
\boldsymbol{\psi }\neq \boldsymbol{\psi }_{0}$ and $\boldsymbol{\psi }\in 
\overline{\Psi }$, we must have $\left( \boldsymbol{\psi }-\boldsymbol{\psi }%
_{0}\right) ^{\prime }\boldsymbol{H}_{N}\left( \boldsymbol{\psi }^{\ast
}\right) \left( \boldsymbol{\psi }-\boldsymbol{\psi }_{0}\right) <0$,
because $\boldsymbol{\psi }\in \overline{\Psi }$ implies $\boldsymbol{\psi }%
^{\ast }\in \overline{\Psi }$ and, therefore, $\boldsymbol{H}\left( 
\boldsymbol{\psi }^{\ast }\right) $ is negative definite. Hence, $L\left( 
\boldsymbol{\psi }\right) <L\left( \boldsymbol{\psi }_{0}\right) $ if $%
\boldsymbol{\psi }\in \overline{\Psi }$ and $\boldsymbol{\psi }\neq 
\boldsymbol{\psi }_{0}$.%

\bigskip

\noindent \textbf{Proof of Theorem 1}: The conclusions of Lemmas B.2 and
B.3 imply there is a measurable maximizer, $\widehat{\boldsymbol{\psi }}$,
in $\overline{\Psi }$ and $\widehat{\boldsymbol{\psi }}\overset{a.s.}{%
\rightarrow }\boldsymbol{\psi }_{0}$ (see, e.g., Amemiya, 1985, Theorem
4.1.1, and his footnote 1 on p. 107).%

\section*{Appendix C: Theorem 2 Proof}

Theorem 2 is proven by establishing several lemmas. The first
result is an elementary inequality, which is applied repeatedly in the
sequel.

\bigskip

\noindent \textbf{Lemma C.1}. For $r>0$, $\left\vert
\sum_{j=1}^{m}a_{j}\right\vert ^{r}\leq b_{r}\sum_{j=1}^{m}\left\vert
a_{j}\right\vert ^{r}$ where $b_{r}=1$ or $2^{\left( r-1\right) \left(
m-1\right) }$ according as $r\leq 1$ or $r\geq 1$.

\bigskip

\noindent \textbf{Proof}. By repeated application of the inequality $%
\left\vert a+b\right\vert ^{r}\leq c_{r}\left\vert a\right\vert
^{r}+c_{r}\left\vert b\right\vert ^{r}$, $r>0$, where $c_{r}=1$ or $2^{r-1}$
according as $r\leq 1$ or $r\geq 1$ (see Lo\`{e}ve 1977, p. 157), we have $%
\left\vert \sum_{j=1}^{m}a_{j}\right\vert ^{r}\leq c_{r}\left\vert
a_{1}\right\vert ^{r}+c_{r}\left\vert \sum_{j=2}^{m}a_{j}\right\vert
^{r}\leq c_{r}\left\vert a_{1}\right\vert ^{r}+c_{r}^{2}\left\vert
a_{2}\right\vert ^{r}+c_{r}^{2}\left\vert \sum_{j=3}^{m}a_{j}\right\vert
^{r}\leq \sum_{j=1}^{m-1}c_{r}^{j}\left\vert a_{j}\right\vert
^{r}+c_{r}^{m-1}\left\vert a_{m}\right\vert ^{r}$. Also, $%
\sum_{j=1}^{m-1}c_{r}^{j}\left\vert a_{j}\right\vert
^{r}+c_{r}^{m-1}\left\vert a_{m}\right\vert ^{r}\leq
b_{r}\sum_{j=1}^{m}\left\vert a_{j}\right\vert ^{r}$ for $b_{r}=c_{r}^{m-1}$.%

\bigskip

\noindent \textbf{Lemma C.2}. Suppose C1$^{\prime }$, C2, C3, C5, and C6 are
satisfied. Then $\sqrt{N}\boldsymbol{g}_{N}\left( \boldsymbol{\psi }%
_{0}\right) \overset{d}{\rightarrow }\mathcal{N}\left( \mathbf{0},\mathcal{I}%
_{0}\right) $.

\bigskip

\noindent \textbf{Proof}. Let $\boldsymbol{\mu }$ be a $m\times 1$ vector of
constants such that $\boldsymbol{\mu }\neq \mathbf{0}$. We have $\boldsymbol{%
\mu }^{\prime }\sqrt{N}\boldsymbol{g}_{N}\left( \boldsymbol{\psi }%
_{0}\right) =N^{-1/2}\sum_{i}\mathcal{Z}_{i}$ for $\mathcal{Z}_{i}=%
\boldsymbol{\mu }^{\prime }\left( \partial l_{i}\left( \boldsymbol{\psi }%
_{0}\right) /\partial \boldsymbol{\psi }\right) $. And $\sqrt{N}\boldsymbol{g%
}_{N}\left( \boldsymbol{\psi }_{0}\right) \overset{d}{\rightarrow }\mathcal{N%
}\left( \mathbf{0},\mathcal{I}_{0}\right) $ if $N^{-1/2}\sum_{i}\mathcal{Z}%
_{i}\overset{d}{\rightarrow }\mathcal{N}\left( 0,\boldsymbol{\mu }^{\prime }%
\mathcal{I}_{0}\boldsymbol{\mu }\right) $ (see Amemiya 1985, Theorem 3.3.8).

\sloppy To verify $N^{-1/2}\sum_{i}\mathcal{Z}_{i}\overset{d}{\rightarrow }%
\mathcal{N}\left( 0,\boldsymbol{\mu }^{\prime }\mathcal{I}_{0}\boldsymbol{%
\mu }\right) $, let $\nu _{i}^{2}=var\left( \mathcal{Z}_{i}\right) =%
\boldsymbol{\mu }^{\prime }E\left[ \left( \partial l_{i}\left( \boldsymbol{%
\psi }_{0}\right) /\partial \boldsymbol{\psi }\right) \left( \partial
l_{i}\left( \boldsymbol{\psi }_{0}\right) /\partial \boldsymbol{\psi }%
\right) ^{\prime }\right] \boldsymbol{\mu }$, and $\overline{\nu }%
_{N}^{2}=N^{-1}\sum_{i}\nu _{i}^{2}$. Because $\lim_{N\rightarrow \infty }%
\overline{\nu }_{N}^{2}=\boldsymbol{\mu }^{\prime }\mathcal{I}_{0}%
\boldsymbol{\mu }$ (by C6), we have $N^{-1/2}\sum_{i}\mathcal{Z}_{i}\overset{%
d}{\rightarrow }\mathcal{N}\left( 0,\boldsymbol{\mu }^{\prime }\mathcal{I}%
_{0}\boldsymbol{\mu }\right) $ if $N^{-1/2}\sum_{i}\mathcal{Z}_{i}/\overline{%
\nu }_{N}\overset{d}{\rightarrow }\mathcal{N}\left( 0,1\right) $. Moreover, $%
N^{-1/2}\sum_{i}\mathcal{Z}_{i}/\overline{\nu }_{N}\overset{d}{\rightarrow }%
\mathcal{N}\left( 0,1\right) $ if $E\left( \mathcal{Z}_{i}\right) =0$, $%
\overline{\nu }_{N}^{2}>\epsilon ^{\prime }>0$ for all $N$ sufficiently
large, and $E\left\vert \mathcal{Z}_{i}\right\vert ^{2+\epsilon /2}<M$ for
all $i$ and some $\epsilon /2>0$ (see White 2001, Theorem 5.10). Therefore,
Lemma C.2 is proven upon proving $E\left( \mathcal{Z}_{i}\right) =0$, $%
\overline{\nu }_{N}^{2}>\epsilon ^{\prime }>0$ for all $N$ sufficiently
large, and $E\left\vert \mathcal{Z}_{i}\right\vert ^{2+\epsilon /2}<M$ for
all $i$ and some $\epsilon /2>0$.

\fussy We can verify $E\left( \mathcal{Z}_{i}\right) =0$ and $\overline{\nu }%
_{N}^{2}>\epsilon ^{\prime }>0$ for all $N$ sufficiently large easily. In
particular, Eq. (\ref{gradient1}) implies $E\left( \mathcal{Z}_{i}\right) =0$%
. Moreover, given C6, we have $\lim_{N\rightarrow \infty }\overline{\nu }%
_{N}^{2}=\boldsymbol{\mu }^{\prime }\mathcal{I}_{0}\boldsymbol{\mu }$, and,
because $\mathcal{I}_{0}$ is positive definite, we can find an $\epsilon
^{\prime }>0$ such that $\overline{\nu }_{N}^{2}>\epsilon ^{\prime }$ for
all $N$ sufficiently large.

To verify $E\left\vert \mathcal{Z}_{i}\right\vert ^{2+\epsilon /2}<M$
for all $i$ and some $\epsilon /2>0$, first let $\mu _{j}$ and $\psi _{j}$
denote the $j$th elements of $\boldsymbol{\mu }$ and $\boldsymbol{\psi }$.
Then $\mathcal{Z}_{i}=\sum_{j}\mu _{j}\partial l_{i}\left( \boldsymbol{\psi }%
_{0}\right) /\partial \boldsymbol{\psi }_{j}$. Hence, by Lemma C.1, we have $%
E\left\vert \mathcal{Z}_{i}\right\vert ^{2+\epsilon /2}<M$ for all $i$ if $%
E\left\vert \partial l_{i}\left( \boldsymbol{\psi }_{0}\right) /\partial 
\boldsymbol{\psi }_{j}\right\vert ^{2+\epsilon /2}<M$ for all $i$ and $j$.
Next, recall $\partial l_{i}\left( \boldsymbol{\psi }_{0}\right) /\partial 
\boldsymbol{\gamma }=\boldsymbol{W}_{i}^{\prime }\Omega _{0}^{-1}\boldsymbol{%
u}_{i}$ while $\partial l_{i}\left( \boldsymbol{\psi }_{0}\right) /\partial 
\boldsymbol{\omega }=-$vech$\left( \Omega _{0}^{-1}-\Omega _{0}^{-1}%
\boldsymbol{u}_{i}\boldsymbol{u}_{i}^{\prime }\Omega _{0}^{-1}\right) /2$.
Moreover, upon letting $\omega _{0}^{st}$ denote the $\left( s,t\right) $th
element of $\Omega _{0}^{-1}$ and recalling $W_{isj}$ denotes the $\left(
s,j\right) $th element of $\boldsymbol{W}_{i}$, the elements of $\boldsymbol{%
W}_{i}^{\prime }\Omega _{0}^{-1}\boldsymbol{u}_{i}$ are of the form $%
\sum_{s}\sum_{t}\omega _{0}^{st}W_{isj}u_{it}$ while the elements of vech$%
\left( \Omega _{0}^{-1}-\Omega _{0}^{-1}\boldsymbol{u}_{i}\boldsymbol{u}%
_{i}^{\prime }\Omega _{0}^{-1}\right) $ are of the form $\omega
_{0}^{jk}-\sum_{s}\sum_{t}\omega _{0}^{js}\omega _{0}^{kt}u_{is}u_{it}$.
These observations and another application of Lemma C.1 implies $E\left\vert
\partial l_{i}\left( \boldsymbol{\psi }_{0}\right) /\partial \boldsymbol{%
\psi }_{j}\right\vert ^{2+\epsilon /2}<M$ for all $i$ and $j$ if $%
E\left\vert W_{isj}u_{it}\right\vert ^{2+\epsilon /2}<M$ and $E\left\vert
u_{is}u_{it}\right\vert ^{2+\epsilon /2}<M$ for all $i$, $j$, $s,$ and $t$.
But $E\left\vert W_{isj}u_{it}\right\vert ^{2+\epsilon /2}\leq \left(
E\left\vert W_{isj}\right\vert ^{4+\epsilon }E\left\vert u_{it}\right\vert
^{4+\epsilon }\right) ^{1/2}$ by the Cauchy-Schwarz inequality. Moreover,
for a suitable choice of $\epsilon >0,$ we have $E\left\vert
W_{isj}\right\vert ^{4+\epsilon }<M$ for all $i,$ $s$, and $j$ by C1$%
^{\prime }$. Condition C1$^{\prime }$ also implies $E\left\vert
u_{it}\right\vert ^{4+\epsilon }<M$ for all $i$ and $t$. Hence, $E\left\vert
W_{isj}u_{it}\right\vert ^{2+\epsilon /2}<M$ for all $i$, $j$, $s,$ and $t$.
Similar arguments give $E\left\vert u_{is}u_{it}\right\vert ^{2+\epsilon
/2}<M$ for all $i$, $s,$ and $t$. It follows that $E\left\vert \mathcal{Z}%
_{i}\right\vert ^{2+\epsilon /2}<M$ for all $i$ and some $\epsilon /2>0$.%

\bigskip

\noindent \textbf{Lemma C.3}. Let $\overline{\Psi }$ be a compact subset of $%
\Psi $. Suppose C1, C4, and C5 are satisfied. Then $\boldsymbol{H}_{N}\left( 
\boldsymbol{\cdot }\right) \overset{a.s.}{\rightarrow }\boldsymbol{H}\left( 
\boldsymbol{\cdot }\right) $ uniformly on $\overline{\Psi }$.

\bigskip

\sloppy \noindent \textbf{Proof}. Let $h_{\gamma _{j}\gamma _{k}}\left( 
\boldsymbol{\psi }\right) =\lim_{N\rightarrow \infty }E\left[ \partial
^{2}L_{N}\left( \boldsymbol{\psi }\right) /\partial \gamma _{j}\partial
\gamma _{k}\right] $. Then $\left\vert \partial ^{2}L_{N}\left( \boldsymbol{%
\psi }\right) /\partial \gamma _{j}\partial \gamma -h_{\gamma _{j}\gamma
_{k}}\left( \boldsymbol{\psi }\right) \right\vert =\left\vert
\sum_{s}\sum_{t}\omega ^{st}\left(
S_{W_{sj}W_{tk},N}+A_{W_{sj}W_{tk},N}\right) \right\vert \leq
\sum_{s}\sum_{t}\left\vert \omega ^{st}\right\vert \left( \left\vert
S_{W_{sj}W_{tk},N}\right\vert +\left\vert A_{W_{sj}W_{tk},N}\right\vert
\right) .$ (For the definitions of $S_{W_{sj}W_{tk},N}$ and $%
A_{W_{sj}W_{tk},N}$, see the proof of Lemma B.2.) Given $\omega ^{st}$ is
bounded for $\boldsymbol{\psi }\in \overline{\Psi }$, we have $\sup_{%
\boldsymbol{\psi }\in \overline{\Psi }}\left\vert \partial ^{2}L_{N}\left( 
\boldsymbol{\psi }\right) /\partial \gamma _{j}\partial \gamma
_{k}-h_{\gamma _{j}\gamma _{k}}\left( \boldsymbol{\psi }\right) \right\vert
\leq M\sum_{s}\sum_{t}\left( \left\vert S_{W_{sj}W_{tk},N}\right\vert
+\left\vert A_{W_{sj}W_{tk},N}\right\vert \right) $. Recall that $%
S_{W_{sj}W_{tk},N}\overset{a.s.}{\rightarrow }0$ (see the proof of Lemma
B.2), and $A_{W_{sj}W_{tk},N}\rightarrow 0$. Therefore, $\partial
^{2}L_{N}\left( \cdot \right) /\partial \gamma _{j}\partial \gamma _{k}%
\overset{a.s.}{\rightarrow }h_{\gamma _{j}\gamma _{k}}\left( \cdot \right) $
uniformly on $\overline{\Psi }$.

Let $h_{\gamma _{j}\omega _{k}}\left( \boldsymbol{\psi }\right)
=\lim_{N\rightarrow \infty }E\left[ \partial ^{2}L_{N}\left( \boldsymbol{%
\psi }\right) /\partial \gamma _{j}\partial \omega _{k}\right] $. Also, let $%
\vartheta _{k,st}$ denote the $\left( s,t\right) $the element of $\Omega
^{-1}\left( \partial \Omega /\partial \omega _{k}\right) \Omega ^{-1}$. Then 
$\partial ^{2}L_{N}\left( \boldsymbol{\psi }\right) /\partial \gamma
_{j}\partial \omega _{k}-h_{\gamma _{j}\omega _{k}}\left( \boldsymbol{\psi }%
\right) =-\sum_{s}\sum_{t}\vartheta
_{k,st}[S_{y_{s}W_{tj},N}+A_{y_{s}W_{tj},N}]+\sum_{s}\sum_{t}\sum_{l}%
\vartheta _{k,st}\gamma _{l}\left[ S_{W_{sj}W_{tl},N}+A_{W_{sj}W_{tl},N}%
\right] $. (For the definitions of $S_{y_{s}W_{tj},N}$ and $%
A_{y_{s}W_{tj},N} $, see the proof of Lemma B.2$.$) Because $\vartheta
_{k,st}$ is a continuous function on $\overline{\Psi }$, and, therefore,
bounded on $\overline{\Psi }$, and $\gamma _{l}$ is bounded for $\boldsymbol{%
\psi }\in \overline{\Psi }$, we have $\sup_{\boldsymbol{\psi }\in \overline{%
\Psi }}\left\vert \partial ^{2}L_{N}\left( \boldsymbol{\psi }\right)
/\partial \gamma _{j}\partial \omega _{k}-h_{\gamma _{j}\omega _{k}}\left( 
\boldsymbol{\psi }\right) \right\vert \leq M\sum_{s}\sum_{t}\left(
\left\vert S_{y_{s}W_{tj},N}\right\vert +\left\vert
A_{y_{s}W_{tj},N}\right\vert \right) +M\sum_{s}\sum_{t}\sum_{l}\left(
\left\vert S_{W_{sj}W_{tl},N}\right\vert +\left\vert
A_{W_{sj}W_{tl},N}\right\vert \right) $. Given $S_{y_{s}W_{tj},N}\overset{%
a.s.}{\rightarrow }0$, $S_{W_{sj}W_{tl},N}\overset{a.s.}{\rightarrow }0$, $%
A_{y_{s}W_{tj},N}\rightarrow 0$, and $A_{W_{sj}W_{tl},N}\rightarrow 0$, we
have $\partial ^{2}L_{N}\left( \cdot \right) /\partial \gamma _{j}\partial
\omega _{k}\overset{a.s.}{\rightarrow }h_{\gamma _{j}\omega _{k}}\left(
\cdot \right) $ uniformly on $\overline{\Psi }$.

Finally, from (\ref{2nd_deriv_3}), we see that $\partial ^{2}L_{N}\left( 
\boldsymbol{\psi }\right) /\partial \omega _{j}\partial \omega _{k}-E\left(
\partial ^{2}L_{N}\left( \boldsymbol{\psi }\right) /\partial \omega
_{j}\partial \omega _{k}\right) =-\left( 2N\right) ^{-1}\sum_{i}\left\{
s_{ijk}^{\left( 1\right) }\left( \boldsymbol{\psi }\right) -E\left[
s_{ijk}^{\left( 1\right) }\left( \boldsymbol{\psi }\right) \right] \right\}
-\left( 2N\right) ^{-1}\sum_{i}\left\{ s_{ijk}^{\left( 2\right) }\left( 
\boldsymbol{\psi }\right) -E\left[ s_{ijk}^{\left( 2\right) }\left( 
\boldsymbol{\psi }\right) \right] \right\} $. Note that%
\begin{equation}
\frac{1}{N}\sum_{i}\left\{ s_{ijk}^{\left( 1\right) }\left( \boldsymbol{\psi 
}\right) -E\left[ s_{ijk}^{\left( 1\right) }\left( \boldsymbol{\psi }\right) %
\right] \right\} =\boldsymbol{S}_{\cdot j}^{\prime }\left( \Omega
^{-1}\otimes \Omega ^{-1}\boldsymbol{U}_{N}\left( \boldsymbol{\gamma }%
\right) \Omega ^{-1}\right) \boldsymbol{S}_{\cdot k}  \label{s_1}
\end{equation}%
where $\boldsymbol{U}_{N}\left( \boldsymbol{\gamma }\right)
=N^{-1}\sum_{i}\left\{ \boldsymbol{u}_{i}\left( \boldsymbol{\gamma }\right) 
\boldsymbol{u}_{i}\left( \boldsymbol{\gamma }\right) ^{\prime }-E\left[ 
\boldsymbol{u}_{i}\left( \boldsymbol{\gamma }\right) \boldsymbol{u}%
_{i}\left( \boldsymbol{\gamma }\right) ^{\prime }\right] \right\} $. Because 
$\boldsymbol{S}_{\cdot j}$ is a vector of zeros and ones, we see that the
right-hand side of (\ref{s_1}) is a sum of the elements of $\Omega
^{-1}\otimes \Omega ^{-1}\boldsymbol{U}_{N}\left( \boldsymbol{\gamma }%
\right) \Omega ^{-1}$. Therefore, if each element of this matrix converges
almost surely to zero uniformly on $\overline{\Psi }$, then $%
N^{-1}\sum_{i}\left\{ s_{ijk}^{\left( 1\right) }\left( \cdot \right) -E\left[
s_{ijk}^{\left( 1\right) }\left( \cdot \right) \right] \right\} \overset{a.s.%
}{\rightarrow }0$ uniformly on $\overline{\Psi }$. Similar arguments can be
used to show $N^{-1}\sum_{i}\left\{ s_{ijk}^{\left( 2\right) }\left( \cdot
\right) -E\left[ s_{ijk}^{\left( 2\right) }\left( \cdot \right) \right]
\right\} \overset{a.s.}{\rightarrow }0$ uniformly on $\overline{\Psi }$.

To see that each element of $\Omega ^{-1}\otimes \Omega ^{-1}\boldsymbol{U}%
_{N}\left( \boldsymbol{\gamma }\right) \Omega ^{-1}$ converges almost surely
to zero uniformly, note that the matrix $\Omega ^{-1}\otimes \Omega ^{-1}%
\boldsymbol{U}_{N}\left( \boldsymbol{\gamma }\right) \Omega ^{-1}$ can be
partitioned into $T\times T$ sub-matrices of the form $\omega ^{lm}\Omega
^{-1}\boldsymbol{U}_{N}\left( \boldsymbol{\gamma }\right) \Omega ^{-1}$ ($%
l=1,\ldots ,T$, $m=1,\ldots ,T$). Furthermore, the $\left( j,k\right) $th
element of $\omega ^{lm}\Omega ^{-1}\boldsymbol{U}_{N}\left( \boldsymbol{%
\gamma }\right) \Omega ^{-1}$ is $\omega ^{lm}\sum_{s}\sum_{t}\omega
^{js}\omega ^{kt}N^{-1}\sum_{i}\left\{ u_{is}\left( \boldsymbol{\gamma }%
\right) u_{it}\left( \boldsymbol{\gamma }\right) -E\left[ u_{is}\left( 
\boldsymbol{\gamma }\right) u_{it}\left( \boldsymbol{\gamma }\right) \right]
\right\} $. And, by familiar arguments, we can show that the absolute value
of this element is no greater than $M\sum_{s}\sum_{t}\left\vert
N^{-1}\sum_{i}\left\{ u_{is}\left( \boldsymbol{\gamma }\right) u_{it}\left( 
\boldsymbol{\gamma }\right) -E\left[ u_{is}\left( \boldsymbol{\gamma }%
\right) u_{it}\left( \boldsymbol{\gamma }\right) \right] \right\}
\right\vert $ for $\boldsymbol{\psi \in }\overline{\Psi }$. Moreover, $%
N^{-1}\sum_{i}\left\{ u_{is}\left( \boldsymbol{\gamma }\right) u_{it}\left( 
\boldsymbol{\gamma }\right) -E\left[ u_{is}\left( \boldsymbol{\gamma }%
\right) u_{it}\left( \boldsymbol{\gamma }\right) \right] \right\}
=S_{y_{s}y_{t},N}-\sum_{q}\gamma
_{q}(S_{y_{s}W_{tq},N}+S_{y_{t}W_{sq},N})+\sum_{q}\sum_{r}\gamma _{q}\gamma
_{r}S_{W_{sq}W_{tr}\,N}$, and, given $\boldsymbol{\gamma }$ is bounded for $%
\boldsymbol{\psi \in }\overline{\Psi }$, we have%
\begin{eqnarray}
&&\sup_{\boldsymbol{\psi \in }\overline{\Psi }}\left\vert \frac{1}{N}%
\sum_{i}\left\{ u_{is}\left( \boldsymbol{\gamma }\right) u_{it}\left( 
\boldsymbol{\gamma }\right) -E\left[ u_{is}\left( \boldsymbol{\gamma }%
\right) u_{it}\left( \boldsymbol{\gamma }\right) \right] \right\} \right\vert
\notag  \label{hess_unif_conv_3} \\
&\leq &\left\vert S_{y_{s}y_{t},N}\right\vert +M\sum_{q}\left( \left\vert
S_{y_{s}W_{tq},N}\right\vert +\left\vert S_{y_{t}W_{sq},N}\right\vert
+\sum_{r}\left\vert S_{W_{sq}W_{tr}\,N}\right\vert \right) .
\label{hess_unif_con_3}
\end{eqnarray}%
Because the right-hand side (\ref{hess_unif_con_3}) $\overset{a.s.}{%
\rightarrow }0$ (see the proof of Lemma B.2), we have $N^{-1}\sum_{i}\left\{
s_{i,jk}^{\left( 1\right) }\left( \cdot \right) -E\left[ s_{i,jk}^{\left(
1\right) }\left( \cdot \right) \right] \right\} \overset{a.s.}{\rightarrow }%
0 $ uniformly on $\overline{\Psi }$. Simliar arguments establish $%
N^{-1}\sum_{i}\left\{ s_{i,jk}^{\left( 2\right) }\left( \cdot \right) -E%
\left[ s_{i,jk}^{\left( 2\right) }\left( \cdot \right) \right] \right\} 
\overset{a.s.}{\rightarrow }0$ uniformly on $\overline{\Psi }$. It follows
that $\partial ^{2}L_{N}\left( \cdot \right) /\partial \omega _{j}\partial
\omega _{k}-E\left[ \partial ^{2}L_{N}\left( \cdot \right) /\partial \omega
_{j}\partial \omega _{k}\right] \overset{a.s.}{\rightarrow }0$ uniformly on $%
\overline{\Psi }$.

Let $h_{\omega _{j}\omega _{k}}\left( \boldsymbol{\psi }\right)
=\lim_{N\rightarrow \infty }E\left[ \partial ^{2}L_{N}\left( \boldsymbol{%
\psi }\right) /\partial \omega _{j}\partial \omega _{k}\right] $. We can
establish $E\left[ \partial ^{2}L_{N}\left( \cdot \right) /\partial \omega
_{j}\partial \omega _{k}\right] \rightarrow h_{\omega _{j}\omega _{k}}\left(
\cdot \right) $ uniformly on $\overline{\Psi }$ by arguments paralleling
those in the last two paragraphs. (For example, in the foregoing
derivations, replace $N^{-1}\sum_{i}E\left[ u_{is}\left( \boldsymbol{\gamma }%
\right) u_{it}\left( \boldsymbol{\gamma }\right) \right] $ with $%
\lim_{N\rightarrow \infty }N^{-1}\sum_{i}E\left[ u_{is}\left( \boldsymbol{%
\gamma }\right) u_{it}\left( \boldsymbol{\gamma }\right) \right] $ and $%
N^{-1}\sum_{i}u_{is}\left( \boldsymbol{\gamma }\right) u_{it}\left( 
\boldsymbol{\gamma }\right) $ with $N^{-1}\sum_{i}E\left[ u_{is}\left( 
\boldsymbol{\gamma }\right) u_{it}\left( \boldsymbol{\gamma }\right) \right] 
$. Also, replace $S_{y_{s}y_{t},N}$, $S_{y_{s}W_{tq},N}$, $S_{y_{t}W_{sq},N}$%
, and $S_{W_{sq}W_{tr}\,N}$ with $A_{y_{s}y_{t},N}$, $A_{y_{s}W_{tq},N}$, $%
A_{y_{t}W_{sq},N}$, and $A_{W_{sq}W_{tr}\,N}$.)

\fussy From the foregoing, we have $\partial ^{2}L_{N}\left( \cdot \right)
/\partial \omega _{j}\partial \omega _{k}\overset{a.s.}{\rightarrow }%
h_{\omega _{j}\omega _{k}}\left( \cdot \right) $ uniformly on $\overline{%
\Psi }$.%

\bigskip

\noindent \textbf{Proof of Theorem 2}: The conclusions of Lemmas C.2 and
C.3, the consistency of $\widehat{\boldsymbol{\psi }}$, the continuity of $%
\boldsymbol{H}\left( \cdot \right) $ at $\boldsymbol{\psi }_{0}$, and the
nonsingularity of $\boldsymbol{H}_{0}=\boldsymbol{H}\left( \boldsymbol{\psi }%
_{0}\right) $ imply $\sqrt{N}\left( \widehat{\boldsymbol{\psi }}-\boldsymbol{%
\psi }_{0}\right) \overset{d}{\rightarrow }\mathcal{N}\left( \mathbf{0},%
\boldsymbol{H}_{0}^{-1}\mathcal{I}_{0}\boldsymbol{H}_{0}^{-1}\right) $ (see
Newey and McFadden 1994, Theorem 3.1).%

\section*{Appendix D: Proof of Theorems 3 and 4}

The proofs of Theorems 3 and 4 are similar to the proofs of
Theorems 1 and 2. For example, Conditions C1 and C2$^{\prime }$ ensure the
linear projection parameters in (\ref{diff_lp}) exist and do not depend on $%
i $ and the errors in $\widetilde{\boldsymbol{u}}_{i}$ are uncorrelated with
the regressors in $\boldsymbol{x}_{i}$. Furthermore, the quasi
log-likelihood $\sum_{i=1}^{N}\widetilde{l}\left( \boldsymbol{\lambda }%
_{0}\right) $ is similar to the quasi log-likelihood $\sum_{i=1}^{N}l\left( 
\boldsymbol{\psi }_{0}\right) $, and, therefore, most of the technical
details are the same as in Appendices B and C and need not be repeated.

However, the conlusions of Theorems 3 and 4 depend on $E\left( \widetilde{%
\boldsymbol{W}}_{i}^{\prime }\Upsilon _{0}^{-1}\widetilde{\boldsymbol{u}}%
_{i}\right) =\mathbf{0}$ being true, and the proof of this result, though
similar to the proof of Lemma 1, differs in some details. Therefore, the
proof of $E\left( \widetilde{\boldsymbol{W}}_{i}^{\prime }\Upsilon _{0}^{-1}%
\widetilde{\boldsymbol{u}}_{i}\right) =\mathbf{0}$ is provided in this
appendix.\bigskip

\noindent \textbf{Lemma D.1}. Suppose $E\left( x_{itk}^{2}\right) <\infty $
and $E\left( y_{it}^{2}\right) <\infty $, for each $i$, $t$, and $k$, and
Conditions C2$^{\prime }$ and C3$^{\prime }$ are satisfied. Then $E\left( 
\widetilde{\boldsymbol{W}}_{i}^{\prime }\Upsilon _{0}^{-1}\widetilde{%
\boldsymbol{u}}_{i}\right) =\mathbf{0}$.

\bigskip

\noindent \textbf{Proof}. Let%
\begin{equation*}
\widetilde{\boldsymbol{Z}}_{i}=\left( 
\begin{array}{cc}
\mathbf{0} & \boldsymbol{I}_{p}\otimes \left( 1,\boldsymbol{x}_{i}^{\prime
}\right) \\ 
\Delta \boldsymbol{X}_{i} & \mathbf{0}%
\end{array}%
\right) .
\end{equation*}%
Given this definition, showing $E\left( \widetilde{\boldsymbol{W}}%
_{i}^{\prime }\Upsilon _{0}^{-1}\widetilde{\boldsymbol{u}}_{i}\right) =%
\mathbf{0}$ consists of showing $E\left( \widetilde{\boldsymbol{Z}}%
_{i}^{\prime }\Upsilon _{0}^{-1}\widetilde{\boldsymbol{u}}_{i}\right) =%
\mathbf{0}$ and $E\left[ \left( \mathbf{0,}\text{ }\Delta \boldsymbol{y}%
_{i,-j}^{\prime }\right) \Upsilon _{0}^{-1}\widetilde{\boldsymbol{u}}_{i}%
\right] =0$ ($j=1,\ldots ,p$). Under the conditions of the lemma, the
elements of $\widetilde{\boldsymbol{Z}}_{i}$ are uncorrelated with the
elements of $\widetilde{\boldsymbol{u}}_{i}$; hence, $E\left( \widetilde{%
\boldsymbol{Z}}_{i}^{\prime }\Upsilon _{0}^{-1}\widetilde{\boldsymbol{u}}%
_{i}\right) =\mathbf{0}$. It remains to show $E\left[ \left( \mathbf{0,}%
\text{ }\Delta \boldsymbol{y}_{i,-j}^{\prime }\right) \Upsilon _{0}^{-1}%
\widetilde{\boldsymbol{u}}_{i}\right] =0$.

This result can be established by arguments similar to those used in the
proof of Lemma 1. Specifically, let $\Delta \boldsymbol{\xi }_{it}=\left(
\Delta y_{it},\Delta y_{i,t-1},\ldots ,\Delta y_{i,t-p+1}\right) ^{\prime },$
$\Delta \boldsymbol{\varsigma }_{it}=\left( \Delta \boldsymbol{x}%
_{it}^{\prime }\boldsymbol{\beta }_{0}+\Delta e_{it},0,\ldots ,0\right)
^{\prime }$ and let $\boldsymbol{F}$ be defined as in (\ref{F_def}). Then we
get $\Delta \boldsymbol{\xi }_{i2}=\boldsymbol{F}\Delta \boldsymbol{\xi }%
_{i1}+\Delta \boldsymbol{\varsigma }_{i2}$; and, for $t>2$, we have $\Delta 
\boldsymbol{\xi }_{it}=\boldsymbol{F}^{t-1}\Delta \boldsymbol{\xi }_{i1}+%
\boldsymbol{F}^{t-2}\Delta \boldsymbol{\varsigma }_{i2}+\cdots +\boldsymbol{F%
}\Delta \boldsymbol{\varsigma }_{i,t-1}+\Delta \boldsymbol{\varsigma }_{it}$%
. Let $f_{rs}^{\left( t\right) }$ denote the ($r,s$)th element of $%
\boldsymbol{F}^{t}$. Then, the preceding implies $\Delta
y_{i2}=f_{11}^{\left( 1\right) }\Delta y_{i1}+f_{12}^{\left( 1\right)
}\Delta y_{i0}+\cdots +f_{1p}^{\left( 1\right) }\Delta y_{i,-p+2}+\Delta 
\boldsymbol{x}_{i2}^{\prime }\boldsymbol{\beta }_{0}+\Delta e_{i2}$; and,
for $t>2$, we have $\Delta y_{it}=f_{11}^{\left( t-1\right) }\Delta
y_{i1}+f_{12}^{\left( t-1\right) }\Delta y_{i0}+\cdots +f_{1p}^{\left(
t-1\right) }\Delta y_{i,-p+2}+f_{11}^{\left( t-2\right) }\left( \Delta 
\boldsymbol{x}_{i2}^{\prime }\boldsymbol{\beta }_{0}+\Delta e_{i2}\right)
+\cdots +f_{11}^{\left( 1\right) }\left( \Delta \boldsymbol{x}%
_{i,t-1}^{\prime }\boldsymbol{\beta }_{0}+\Delta e_{i,t-1}\right) +\Delta 
\boldsymbol{x}_{it}^{\prime }\boldsymbol{\beta }_{0}+\Delta e_{it}$ (see the
proof of Lemma 1).

Using these equations we can write $\Delta \boldsymbol{y}_{i,-j}$ as $\Delta 
\boldsymbol{y}_{i,-j}=\widetilde{\boldsymbol{A}}_{j}\Delta \boldsymbol{\xi }%
_{i1}+\widetilde{\boldsymbol{B}}_{j}\left( \Delta \boldsymbol{X}_{i}%
\boldsymbol{\beta }_{0}+\Delta \boldsymbol{e}_{i}\right) $, where $%
\widetilde{\boldsymbol{A}}_{j}$ is a $\left( T-1\right) \times p$ matrix
consisting of the first $T-1$ rows of $\boldsymbol{A}_{j}$ (see Eq. (\ref%
{A_j})) and $\widetilde{\boldsymbol{B}}_{j}$ is a $\left( T-1\right) \times
\left( T-1\right) $ matrix consisting of the first $T-1$ rows and first $T-1$
columns of $\boldsymbol{B}_{j}$ (see Eq. (\ref{B_j})). Recall $\left( \Delta
y_{i,-p+2},\ldots ,\Delta y_{i1}\right) ^{\prime }=\left[ \boldsymbol{I}%
_{p}\otimes \left( 1,\boldsymbol{x}_{i}^{\prime }\right) \right] \boldsymbol{%
\pi }_{0}+\boldsymbol{r}_{i}$ for $\boldsymbol{\pi }_{0}=\left( \mu _{01},%
\boldsymbol{\theta }_{01}^{\prime },\mu _{02},\boldsymbol{\theta }%
_{02}^{\prime },\ldots ,\mu _{0,p},\boldsymbol{\theta }_{0,p}^{\prime
}\right) ^{\prime }$ and $\boldsymbol{r}_{i}=\left( r_{i1},\ldots
,r_{ip}\right) ^{\prime }$ (see Eq. (\ref{levels_sys})). Moreover, note that 
$\Delta \boldsymbol{\xi }_{i1}=\boldsymbol{I}^{\ast }\left( \Delta
y_{i,-p+2},\ldots ,\Delta y_{i1}\right) ^{\prime }$ for $p\times p$ matrix%
\begin{equation*}
\boldsymbol{I}^{\ast }=\left( 
\begin{array}{cccc}
0 & \cdots & 0 & 1 \\ 
0 & \cdots & 1 & 0 \\ 
\vdots &  & \vdots & \vdots \\ 
1 & \cdots & 0 & 0%
\end{array}%
\right) .
\end{equation*}%
Let 
\begin{equation*}
\boldsymbol{D}_{j}=\left( 
\begin{array}{cc}
\mathbf{0} & \mathbf{0} \\ 
\widetilde{\boldsymbol{A}}_{j}\boldsymbol{I}^{\ast } & \widetilde{%
\boldsymbol{B}}_{j}%
\end{array}%
\right) .
\end{equation*}%
Then some straightforward calculations give%
\begin{equation}
\left( 
\begin{array}{cc}
\mathbf{0}, & \Delta \boldsymbol{y}_{i,-j}^{\prime }%
\end{array}%
\right) \Upsilon _{0}^{-1}\widetilde{\boldsymbol{u}}_{i}=\left( \boldsymbol{%
\beta }_{0}^{\prime },\boldsymbol{\pi }_{0}^{\prime }\right) \widetilde{%
\boldsymbol{Z}}_{i}^{\prime }\boldsymbol{D}_{j}^{\prime }\Upsilon _{0}^{-1}%
\widetilde{\boldsymbol{u}}_{i}+\widetilde{\boldsymbol{u}}_{i}^{\prime }%
\boldsymbol{D}_{j}^{\prime }\Upsilon _{0}^{-1}\widetilde{\boldsymbol{u}}_{i}.
\label{y_(i,-j)}
\end{equation}

\sloppy Because the elements of $\widetilde{\boldsymbol{u}}_{i}$ are uncorrelated
with the elements of $\widetilde{\boldsymbol{Z}}_{i}$, we have $E\left[
\left( \boldsymbol{\beta }_{0}^{\prime },\boldsymbol{\pi }_{0}^{\prime
}\right) \widetilde{\boldsymbol{Z}}_{i}^{\prime }\boldsymbol{D}_{j}^{\prime
}\Upsilon _{0}^{-1}\widetilde{\boldsymbol{u}}_{i}\right] =0$. Also, $E\left( 
\widetilde{\boldsymbol{u}}_{i}^{\prime }\boldsymbol{D}_{j}^{\prime }\Upsilon
_{0}^{-1}\widetilde{\boldsymbol{u}}_{i}\right) =$ tr$\left[ \Upsilon
_{0}^{-1}E\left( \widetilde{\boldsymbol{u}}_{i}\widetilde{\boldsymbol{u}}%
_{i}^{\prime }\right) \boldsymbol{D}_{j}^{\prime }\right] =$ tr$\left( 
\boldsymbol{D}_{j}^{\prime }\right) $. But tr$\left( \boldsymbol{D}%
_{j}^{\prime }\right) =0$, because the upper left-hand submatrix $\mathbf{0}$
in $\boldsymbol{D}_{j}$ is square with zeros down its main diagonal and $%
\widetilde{\boldsymbol{B}}_{j}$ is a square matrix with zeros down its main
diagonal, and, therefore, $\boldsymbol{D}_{j}$ has zeros down its main
diagonal. These observations and Eq. (\ref{y_(i,-j)}) prove $E\left[ \left( 
\begin{array}{cc}
\mathbf{0}, & \Delta \boldsymbol{y}_{i,-j}^{\prime }%
\end{array}%
\right) \Upsilon _{0}^{-1}\widetilde{\boldsymbol{u}}_{i}\right] =0$.%

\newpage
\fussy
\noindent \textbf{\Large References:}

\begin{description}
\item Alvarez, J., Arellano, M. (2004). Robust likelihood estimation of
dynamic panel data models. CEMFI Working Paper 0421.%

\item Amemiya, T. (1985). \textit{Advanced Econometrics}. Cambridge, MA:
Harvard University Press.

\item Anderson, T. W., Hsiao, C. (1981). Estimation of dynamic models with
error components. \textit{Journal of the American Statistical Association}
76, 598-606.

\item Arellano, M., Bond, S. (1991). Some tests of specification for panel
data: Monte Carlo evidence and an application to employment equations.\ 
\textit{The Review of Economic Studies} 58, 277-297.%

\item Binder, M., Hsiao, C., Pesaran, M. H. (2005). Estimation and inference
in short panel vector autoregressions with unit roots and cointegration.
Econometric Theory 21, 795-837.

\item Blundell, R., Bond, S. (1998). Initial conditions and moment
restrictions in dynamic panel data models.\ \textit{Journal of Econometrics}
87, 115-143.

\item Bun, M. J. G., Windmeijer, F. (2010). The weak instrument problem of
the system GMM estimator in dynamic panel data models. \textit{The
Econometrics Journal} 13, 95-126.%

\item Chamberlain, G. (1982). Multivariate regression models for panel data. 
\textit{Journal of Econometrics} 18, 5-46.%

\item Chamberlain, G. (1984). Panel data. In: Griliches, Z., Intriligator,
M. D. (eds.), \textit{Handbook of Econometrics}, Vol. 2. Amsterdam: North
Holland, pp. 1247--1318.

\item Hamilton, J. D. (1994). \textit{Time Series Analysis}. Princeton, NJ:
Princeton University Press.

\item Hsiao, C., Pesaran, H. M., Tahmiscioglu, A. K. (2002). Maximum
likelihood estimation of fixed effects dynamic panel data models covering
short time periods. \textit{Journal of Econometrics} 109, 107-150.%

\item Kruiniger, H. (2013). Quasi ML estimation of the panel AR(1)\ model
with arbitrary initial conditions.\ \textit{Journal of Econometrics} 173,
175-188.

\item Lo\`{e}ve, M. (1977). \textit{Probability Theory I}, 4th ed. New York,
NY: Springer-Verlag.

\item Meng, X.-L., van Dyk, D. (1998). Fast EM-Type Implementations for
Mixed Effects Models. \textit{Journal of the Royal Statistical Society,
Series B (Statistical Methodology)} 60, 559-578.

\item Moral-Benito, E. (2013). Likelihood-based estimation of dynamic panels
with predetermined regressors. \textit{Journal of Business \& Economic
Statistics} 31, 451-472.

\item Newey, W. K., McFadden, D. (1994). Large sample estimation and
hypothesis testing.\ In: Engle, R. F., McFadden, D. L. (eds.), \textit{%
Handbook of Econometrics}, Vol. 4. Amsterdam: North Holland, pp. 2111-2245.%

\item Phillips, R. F. (2004). Estimation of a generalized random-effects
model: Some ECME algorithms and Monte Carlo evidence. \textit{Journal of
Economic Dynamics \& Control} 28, 1801-1824.%

\item Phillips, R. F. (2010). Iterated feasible generalized least-squares
estimation of augmented dynamic panel data models.\ \textit{Journal of
Business \& Economic Statistics} 28, 410-422.%

\item Phillips, R. F. (2012). On computing maximum-likelihood estimates of
the unbalanced two-way random-effects model. \textit{Communications in
Statistics--Simulation and Computation} 41, 1921-1927.%

\item Phillips, R. F. (2015). On quasi maximum-likelihood estimation of
dynamic panel data models. \textit{Economics Letters} 137, 91-94.%

\item Rao, C. R. (1973). \textit{Linear Statistical Inference and its
Applications}. New York, NY: Wiley \& Sons.%

\item Roodman, D. (2009). A note on the theme of too many instruments. 
\textit{Oxford Bulletin of Economics and Statistics} 71, 135-158.%
\

\item Ruud, P. A. (2000). \textit{An Introduction to Classical Econometric
Theory}. New York, NY: Oxford University Press.%

\item White, H. (2001). \textit{Asymptotic Theory for Econometricians}. New
York, NY: Academic Press.

\item Wooldridge, J. M. (2010). \textit{Econometric Analysis of Cross
Section and Panel Data}, 2nd ed. Cambridge, MA: MIT Press.%

\end{description}

\end{document}